\documentclass[12pt,reqno]{amsart}
\usepackage{amssymb,delarray}
\usepackage{amsfonts,amscd}
\usepackage{epsfig}
\usepackage[all]{xy}

\def\leq{\leqslant}
\def\geq{\geqslant}

\newcommand{\Sing}{\operatorname{Sing}}

\newtheorem{thm}{Theorem}
{Lemma}
\newtheorem{prop}
{Proposition}
{Problem}
\newtheorem{claim}
{Claim}
\newtheorem{df}
{Definition}
\newtheorem{cor}
{Corollary}
\newtheorem{rem}
{Remark}
{Question}
\newtheorem{ex-thm}{Theorem-Example}

{\catcode`\@=11
\gdef\n@te#1#2{\leavevmode\vadjust{%
 {\setbox\z@\hbox to\z@{\strut#1}%
  \setbox\z@\hbox{\raise\dp\strutbox\box\z@}\ht\z@=\z@\dp\z@=\z@%
  #2\box\z@}}}
\gdef\leftnote#1{\n@te{\hss#1\quad}{}}
\gdef\rightnote#1{\n@te{\quad\kern-\leftskip#1\hss}{\moveright\hsize}}
\gdef\?{\FN@\qumark}
\gdef\qumark{\ifx\next"\DN@"##1"{\leftnote{\rm##1}}\else
 \DN@{\leftnote{\rm??}}\fi{\rm??}\next@}}

\begin{document}
\baselineskip=13.7pt plus 2pt 

\title[On rigid
germs of finite morphisms of smooth surfaces.
]{On rigid
germs of finite morphisms of smooth surfaces. 
}

\author[Vik.S. Kulikov]{Vik.S. Kulikov}

\address{Steklov Mathematical Institute of Russian Academy of Sciences, Moscow, Russia}
 \email{kulikov@mi.ras.ru}

\dedicatory{} \subjclass{}
\thanks{ 
This work is supported by the Russian Science Foundation under grant no.  19-11-00237.}


\maketitle

\def\st{{\sf st}}

\quad \qquad \qquad

\begin{abstract} In the article, we establish and investigate the correspondence between the set of rigid germs of finite morphisms  branched in germs of curves having  $ADE$ singularity types and the set of Belyi rational functions  $f\in \overline{\mathbb Q}(z)$.
\end{abstract} \vspace{0.5cm}

\def\st{{\sf st}}

\setcounter{section}{-1}

\section{Introduction}
Before to formulate the main results of this article, we briefly recall the (well-known) definitions and facts related to germs of  finite morphisms of smooth surfaces.

Let $h_1(u,v)$ and $h_2(u,v)$ be holomorphic functions  in  open subsets $V_i\subset \mathbb C^2$, $i=1,2$. We will assume that $o=(0,0)\in V_1\cap V_2$ and $h_1(0,0)=h_2(0,0)=0$. The functions $h_1(u,v)$ and $h_2(u,v)$ are said to be {\it equivalent} if $h_1(u,v)=h_2(u,v)$ in $V_1\cap V_2$ and their equivalence class is called the {\it germ}  of the function $h(u,v):=h_i(u,v)$ at the point $o$. We can choose $\varepsilon_1\ll 1$, $\varepsilon_2\ll 1$ such that  the closure $\overline{\mathbb D}^2_{\varepsilon_1,\varepsilon_2}$ of the bi-disk $\mathbb D^2_{\varepsilon_1,\varepsilon_2}=\{ (u,v)\in\mathbb C^2 \mid |u|<\varepsilon_1,\,\, |v|<\varepsilon_2 \}$ is contained in $V_1\cap V_2$  and call $h(u,v)$, defined in
$\mathbb D^2_{\varepsilon_1,\varepsilon_2}$
a {\it representative of the germ} (or simply, a {\it germ}) of holomorphic function.
Note that the germ $h(u,v)$ can be given  as a power series
$$ h(u,v)=\sum_{i=0}^{\infty}\sum_{j=0}^{\infty}a_{i,j}u^iv^j\in \mathbb C[[u,v]]$$
absolutely converging in
$\overline{\mathbb D}^2_{\varepsilon_1,\varepsilon_2}$.

A germ of divisor given in $\mathbb D^2_{\varepsilon_1,\varepsilon_2}$ by equation $h(u,v)=0$ is called a {\it curve germ} if it has not multiple components.
Let $B_1,\dots , B_k$ be the irreducible components of a curve germ $B\subset V$
and let $\sigma: V_n\to V$, $\sigma=\sigma_1\circ \dots \circ \sigma_n$, be the minimal sequence of $\sigma$-processes $\sigma_i:V_i\to V_{i-1}$ with centers at points resolving the singular point $o$ of the germ $B$ and such that $\sigma^{-1}(B)$ is a divisor with normal crossings. Denote by $E_{k+i}\subset V_n$ the proper inverse image of the exceptional curve of the blowup
$\sigma_i$ and by $B_j'\subset V_n$ the proper inverse image of the germ $B_j$.

\begin{df} \label{grph} The {\it graph $\Gamma(B)$ of a curve germ} $B$ is a weighted graph having $n+k$ vertices $v_i$. Its vertices $v_i:=b_i$, $i=1,\dots, k$ are in one-to-one correspondence with the curve germs $B'_1,\dots, B'_k$ and their weights are $w_i=0$; the vertices $v_{i+k}:=e_{i+k}$, $i=1,\dots, n$,  are in one-to-one correspondence with the curves   $E_{1+k},\dots, E_{n+k}$ and their weights are $w_{i+k}=(E_{i+k}^2)_{V_n}$; vertices $v_i$ and $v_j$ are connected by an edge in $\Gamma(B)$ if and only if the corresponding to them curves and curve germs have non-empty intersection.
\end{df}

By definition, a {\it family} $h_{\tau}(u,v)$ of {\it germs of functions}, {\it parameterized by the points of a closed disk}   $\overline{\mathbb D}_{\delta}=\{ \tau\in\mathbb C\mid |\tau|\leq \delta\}$,   is a holomorphic function
\begin{equation}\label{deform} h_{\tau}(u,v):=h(u,v,\tau)=\sum_{n=1}^{\infty}\,\, \sum_{i+j=n}a_{i,j}(\tau)u^iv^{j}\end{equation}
defined in $\mathbb D^2_{\varepsilon_1,\varepsilon_2}\times \overline{\mathbb D}_{\delta}$ such that $(\mathcal V=\mathbb D^2_{\varepsilon_1,\varepsilon_2}\times
\overline{\mathbb D}_{\delta}, \mathcal B, \text{pr}_2)$ is a {\it family of curve germs}, where $\mathcal B=(h_{\tau}(u,v))_{\text{red}}$ is the reduced divisor in $\mathcal V$ of the function $h(u,v,\tau)$  and the  restriction to $\mathcal B$ of the projection $\text{pr}_2:\mathcal V\to \overline{\mathbb D}_{\delta}$ is a flat holomorphic map.

\begin{df} \label{def1} {\rm (\cite{Z1}, see also \cite{W})} A  family $h_{\tau}(u,v)$ of germs of functions, parameterized by the points of a closed disk
$\overline{\mathbb D}_{\delta}$ {\rm (}resp., the family $(\mathcal V,\mathcal B,\text{pr}_2)$ given by the family $h_{\tau}(u,v)${\rm )}, is a {\it strong equisingular deformation} if the family $(\mathcal V, \mathcal B, \text{pr}_2)$ of its divisors is a {\it strong equisingular deformation} of curve germs, that is $\Sing\, \mathcal B=\{ o\}\times \overline{\mathbb D}_{\delta}$ and there exists a finite sequence of  
monoidal transformations {\rm (}blowups{\rm )}  $\widetilde{\sigma}_i:\mathcal V_{i}\to \mathcal V_{i-1}$, $i=1,\dots,n$ {\rm (}where $\mathcal V_0=\mathcal V${\rm )},  with centers in smooth curves $\mathcal S_i\subset \Sing \,  \mathcal B_i$, where  $\mathcal B_0=\mathcal B$ and $\mathcal B_{i+1}=\widetilde{\sigma}_{i}^{-1}(\mathcal B_i)$, and such that
\begin{itemize}
\item[$(i)$] $Sing\, \mathcal B_i$ is a disjoint union of sections of $\text{pr}_2\circ\widetilde{\sigma}_1\circ\dots\circ\widetilde{\sigma}_{i-1}$  for each $i$,
\item[$(ii)$] $\mathcal B_{n}$ is a divisor with normal crossings in  $\mathcal V_{n}$.
\end{itemize}
\end{df}

We say that a strong equisingular deformation $h_{\tau}(u,v)$ is {\it trivial} if $h_{\tau}(u,v)$ does not depend on $\tau$.

It is said that two  germs of functions $g_1(u,v)$ and $g_2(u,v)$ (resp. the germs of their divisors)  have the same {\it singularity type} if there is a finite sequence $g_1(u,v)=h_1(u,v),\dots, h_n(u,v)=g_2(u,v)$ of germs of functions such that $h_i(u,v)$ and $h_{i+1}(u,v)$ are members of strong equisingular deformations for $i=1,\dots, n-1$.
Denote by  $T[h(u,v)]$ the singularity type of the  germ of function  $h(u,v)$. We have the following
\begin{prop} \label{Zar} {\rm (}\cite{Z1}, \cite{W}{\rm )}  Two curve germs $(B_1,o)$ and $(B_2,o)$ have the same singularity type if and only if the graphs $\Gamma(B_1)$ and $\Gamma(B_2)$ are isomorphic as weighted graphs.
\end{prop}

\begin{df} \label{rigid} A  germ  of function $h(u,v)$ is rigid if for each germ of function $h_1(u,v)$ such that  $T[h_1(u,v)]=T[h(u,v)]$ there exists a coordinate change  $u_1=u_1(u,v),\,\, v_1=v_1(u,v)$ in $V$ such that $h(u,v)=h_1(u_1(u,v),v_1(u,v))$.
\end{df}

\begin{prop} {\rm (\cite{Ar})} \label{Anrig} A germ of function $h(u,v)$ is rigid if and only if it has one of the following $ADE$ singularity types:

\begin{itemize}
\item $A_n:=T[u^2-v^{n+1}
]$, $n\geq 0$;
\item $D_n:=T[v(u^2-v^{n-2})
]$, $n\geq 4$;
\item $E_6:=T[u^3-v^4
]$;
\item $E_7:=T[u(u^2-v^3)
]$;
\item $E_8:=T[u^3-v^5
]$.
\end{itemize}
\end{prop}

Let $(V,o)=(\mathbb D^2_{\varepsilon_1,\varepsilon_2},o)=\{ (u,v)\in \mathbb C^2\mid |u|<\varepsilon_1,\,\, |v|<\varepsilon_2\}$ be a bi-disk in $\mathbb C^2$, $(U,o')$ a connected germ of smooth complex-analytic surface, and $F: (U,o')\to (V,o)$ a germ of finite holomorphic mapping (below,  a {\it germ of cover}) of local degree $\deg_{o'} F=d$, given in local coordinates $z,w$ in $(U,o')$ by two representatives of germs of holomorphic in $U$ functions
$$  u= f_1(z,w), \qquad  v= f_2(z,w). $$
Denote by $R\subset (U,o')$ the germ of ramification divisor of $F$ given by equation
$$J(F):= \det \left(\begin{array}{cc} \frac{\partial u}{\partial z} & \frac{\partial u}{\partial w} \\ \frac{\partial v}{\partial z} & \frac{\partial v}{\partial w}\end{array}\right)= 0 $$
and by $B=F(R_{red})\subset (V,o)$ the germ of branch curve of $F$. Note that the germ $R\subset (U,o')$ and the curve germ $B\subset (V,o)$ depend only on $F$ and do not depend on the choice of coordinates in $(U,o')$ and $(V,o)$.

We say that two germs $F_1: (U_1,o')\to (V,o)$ and $F_2: (U_2,o')\to (V,o)$ of finite covers are {\it equivalent} if there is a neighbourhood $(W,o)\subset (V,o)$ and bi-holomorphic mappings $G_1:W\to W$ and $G_2:\widetilde W_1\to\widetilde W_2$ such that the following diagram $$ \begin{CD} \widetilde W_1 @>G_2>>  
\widetilde W_2 \\
@VF_1VV  @VVF_2V \\
W  @ >>G_1>  W \end{CD} $$
is commutative, where $\widetilde W_i=F_i^{-1}(W)$.

\begin{df} \label{def-d} We say that a finite holomorphic mapping $\mathcal F: \mathcal U\to (V,o)\times \overline{\mathbb D}_{\delta}$ from a smooth three dimensional complex manifold $\mathcal U$, $\deg \mathcal F =d$, branched along a surface $\mathcal B\subset (V,o)\times \overline{\mathbb D}_{\delta}$,
is a strong deformation of a germ of cover $F_0=\mathcal F_{\mid (U_0,o'_0)}: (U_0,o'_0)\to (V,o)\times \{ \tau=0\}$ and the germs of covers  $F_{\tau_0}=\mathcal F_{\mid (U_{\tau_0},o'_{\tau_0})}: (U_{\tau_0},o'_{\tau_0})\to (V,o)\times \{ \tau=\tau_0 \}$ are strong deformation equivalent to the germ $F_0$, where $U_{\tau_0}=\mathcal F^{-1}(V\times\{\tau=\tau_0\}$, if
\begin{itemize}
\item[$(i)$] the differential gorm $\mathcal F^*(d\tau)\neq 0$ in each point $p\in\mathcal U$,
\item[$(ii)$] $o'_{\tau_0}=\mathcal F^{-1}(o\times\{\tau=\tau_0\})$ is a point for each $\tau_0\in\overline{\mathbb D}_{\delta}$;
\item[$(iii)$]  
$((V,o)\times \overline{\mathbb D}_{\delta}, \mathcal B, \text{pr}_2)$ is a strong equisingular deformation of the curve germ
$B_0=\mathcal B\cap \text{pr}_2^{-1}(0)$.
\end{itemize}
\end{df}

It is said that two  germs of covers  $G_1: (U_1,o'_1)\to (V,o)$ and $G_2: (\widetilde U,\widetilde o')\to (V,o)$ are {\it deformation equivalent}  if there is a finite sequence $$G_1=F_1:(U_1,o'_1)\to (V,o),\dots, F_n=G_2:(\widetilde U,\widetilde o')\to (V,o)$$ of germs of finite covers such that $F_i$ and $F_{i+1}$ are strongly deformation equivalent 
for $i=1,\dots, n-1$.

\begin{df} A germ of cover $F:(U,o')\to (V,o)$ is {\it rigid} if any  deformation equivalent to $F$ cover $F_1:(U_1,o'_1)\to (V,o)$ is equivalent to it.
\end{df}

In Section 1, we prove the following
\begin{thm} \label{main-t} If the germ $B$ of the branch curve of a finite cover $F:(U,o')\to (V,o)$ has one of $ADE$-singularity types, then  $F$ is a rigid germ of finite cover.
\end{thm}

A germ of cover $F$ of degree $\deg F=d$ defines a homomorphism $F_*:\pi_1(V\setminus B,p)\to \mathbb S_d$ (the {\it monodromy} of the germ $F$), where $\mathbb S_d$ is the symmetric group acting on the fibre $F^{-1}(p)$.  The group $G_F=\text{im} F_*$ is called the ({\it local}) {\it monodromy group} of $F$. Note that $G_F$ is a transitive subgroup of $\mathbb S_d$.

Denote by $\mathcal R=(\bigcup_{n\geq 1}\mathcal R_{A_n})\cup (\bigcup_{n\geq 4}\mathcal R_{D_n})\cup(\bigcup_{n\in\{6,7,8\}}\mathcal R_{E_n})$ the set of rigid germs of finite covers brunched along curve germs having, resp., the singularity types $A_n$, $n\geq 1$, $D_n$, $n\geq 4$,  and $E_6$, $E_7$, $E_8$. It follows from Theorem 6 in \cite{K} that $\mathcal R_{T}\neq \emptyset $ for each $ADE$-singularity type $T$. In Section 3, a connection of the set $\mathcal R$  with the set $\mathcal Bel$ of rational Belyi functions is investigated. A rational function $f:\mathbb P^1\to\mathbb P^1$, defined over the algebraic closure $\overline{\mathbb Q}$ of the field of rational function $\mathbb Q$,
is called a {\it Belyi function} if it has no more than three critical values, $\mathcal Bel=\mathcal Bel_{2}\cup \mathcal Bel_3$, where  $\mathcal Bel_{2}$ is the set of Belyi functions with no more than two critical values and the Belyi functions $f\in \mathcal Bel_3$ have three critical values. Further, without loss of generality (due to the action of the group $\text{PGL}(2,\mathbb C)$ on $\mathbb P^1$), we will assume that  $f\in \mathcal Bel_{2}$ are functions  $z=x^n$ and the set of their critical values  is $B_f=\{ 0,\infty\}$ (if $n\geq 1$) and  the set  of critical values of $f\in \mathcal Bel_{3}$ is
$B_f=\{ 0,1,\infty\}$.

Similar to the two-dimensional case, a function $f\in\mathcal Bel$ defines a homomorphism $f_*:\pi_1(\mathbb P^1\setminus B_f,p
) \to \mathbb S_n$ (the {\it monodromy} of $f$), where $n=\deg f$. The image $G_f=\text{im} f_*$ is called the  {\it monodromy group} of $f$. If $f\in\mathcal Bel_{2}$, then $G_f= \mathbb Z_n\subset\mathbb S_n$ is a cyclic group of order $n$.

The group $\pi_1(\mathbb P^1\setminus \{ 0,1,\infty\},p)$ is the free group generated by two 
loops $\gamma_0$ and $\gamma_1$ around the points $0$ and $1$ such that the loop $\gamma_{\infty}=\gamma_0\gamma_1$ is the trivial element in $\pi_1(\mathbb P^1\setminus \{ 0,1\},p)$.
For $f\in \mathcal Bel_3$ denote by $$T_c(f)=\{ c_i=(m_{1,i},\dots ,m_{k_i,i})\}_{m_{1,i}+\dots + m_{k_i,i}=\deg f,\,\, i\in \{ 0,1,\infty\}}$$
the the set of cycle types of permutations $f_*(\gamma_i)$. 
Then, by Hurwitz formula connecting the degree $\deg f$ of $f:\mathbb P^1\to \mathbb P^1$ and the orders of ramification at the critical points of $F$,
\begin{equation} \label{H} n+2=k_0+k_1+k_{\infty}.\end{equation}
Conversely, if a transitive group $G\subset \mathbb S_n$ is generated by two permutations $\sigma_0$ and $\sigma_1$ such that their cycle types and the cycle type of $\sigma_{\infty}=\sigma_0\sigma_1$ satisfy equation (\ref{H}) then there is a rational Belyi function $f$ such that $f_*(\gamma_i)=\sigma_i$.

In Section 2, the Mumford presentations of local fundamental groups $\pi_1^{loc}(B,o)$ is given for the curve germs $B$ having $ADE$ singularity types and which are used in Section 3 in order to prove
\begin{thm} \label{main2} There exists a natural map $\beta:\mathcal R\to \mathcal Bel$ {\rm (see Definition \ref{Be} in Section 3)} such that
\begin{itemize} \item[($\text{\rm i}$)] $\deg \beta(F)$ is a divisor of $\deg_{o'} F$,
\item[($\text{\rm ii}$)] 
$\beta$ induces epimorphisms $\beta_*: G_F\to G_{\beta(F)}$,
\item[($\text{\rm iii}$)]   $\beta(\mathcal R_{A_0}\cup \mathcal R_{A_1})=\mathcal Bel_{2}$,
\item[($\text{\rm iv}$)]  
for each $f\in\mathcal Bel_3$ of degree $\deg f=n$, there exists a germ of cover $F\in \mathcal R_{D_4}\cap \beta^{-1}(f)$ of degree $\deg_{o'} F=n^2$.\end{itemize}
\end{thm}

In Section 4, we prove the following
\begin{thm} \label{main3} For each germ of finite cover $F\in \mathcal R$, there exists an equivalent to $F$ germ $F_1:(U_1,o'_1)\to (V,o)$ given by
$$ u=\frac{f_1(z,w)}{g_1(z,w)},\quad v=\frac{f_2(z,w)}{g_2(z,w)},$$
where $f_i(z,w)$ and $g_i(z,w)\in \overline{\mathbb Q}[z,w]$ for $i=1,2$ and $g_1(0,0)g_2(0,0)\neq 0$.
\end{thm}

\section{Proof of Theorem \ref{main-t}}
Below, the phrase "$f(\tau)$ is a function" means that $f(\tau)$ is a holomorphic function in some neighborhood of the closure  $\overline{\mathbb D}_{\delta}$ of a disk
$\mathbb D_{\delta}=\{ \tau\in\mathbb C\mid |\tau|<\delta\}$.

If a function $f(\tau)$ is not identically equal to zero, then there are at most finitely many points $\alpha_1,\dots, \alpha_n\in \overline{\mathbb D}_{\delta}$ at which $f(\tau)$ vanishes. Let $k_i$ be the order of zero of $f(\tau)$ at $\alpha_i$. Then $f(\tau)=h_f(\tau)p_f(\tau)$, where $p_f(\tau)= \prod_{i=1}^{n}(\tau-\alpha_i)^{k_i}$ and $h_f(\tau)$ is a function invertible in $\overline{\mathbb D}_{\delta}$.  The functions $p_f(\tau)$, $h_f(\tau)$ and the presentation $f(\tau)=h_f(\tau)p_f(\tau)$ will be called respectively the {\it polynomial} and {\it invertible parts}, and {\it canonical factorisation} of $f(\tau)$ in $\mathbb D_{\delta}$.   For functions $f(\tau)$ and $g(\tau)$ such that $f(\tau)g(\tau)\not\equiv 0$, we put $$\text{G.C.D.}(f(\tau),g(\tau))_{\delta}:=\text{G.C.D.}(p_f(\tau),p_g(\tau)),$$ where $p_f(\tau)$ and $p_g(\tau)$ are the polynomial parts of $f(\tau)$ and $g(\tau)$ in $\mathbb D_{\delta}$.

\begin{claim} \label{line} Let $f_1(\tau)$ and $f_2(\tau)$ be two functions such that $\text{G.C.D.}(f_1(\tau),f_2(\tau))_{\delta}=1$.
Then there are functions $g_1(\tau)$ and $g_2(\tau)$
such that $f_1(\tau)g_1(\tau)+f_2(\tau)g_2(\tau)\equiv 1$.
\end{claim}
\proof Let $f_i(\tau)=h_{f_i}(\tau)p_{f_i}(\tau), i=1,2$, be the canonical factorizations in $\mathbb D_{\delta}$. Then
$\text{G.C.D.}(p_{f_1}(\tau),p_{f_2}(\tau))=1$.
Therefore there are polynomials $q_1(\tau)$ and $q_2(\tau)$ such that
$p_{f_1}(\tau)q_1(\tau)+p_{f_2}(\tau)q_2(\tau)\equiv 1$, since $\mathbb C[\tau]$ is a principal ideal ring.
Hence, $g_1(\tau)=\frac{q_1(\tau)}{h_{f_1}(\tau)}$ and $g_2(\tau)=\frac{q_2(\tau)}{h_{f_2}(\tau)}$ are desired functions. \qed \\

Let $(u,v)$ be homogeneous coordinates in $\mathbb P^1$,  $s: \overline{\mathbb D}_{\delta}\to \mathbb P^1\times\mathbb C$ a holomorphic map from a disk $\overline{\mathbb D}_{\delta}\subset \mathbb C$ such that  $\text{pr}_2\circ s=\text{id}: \overline{\mathbb D}_{\delta}\to \overline{\mathbb D}_{\delta}$,  $S=s(\overline{\mathbb D}_{\delta})$ the   section over
$\overline{\mathbb D}_{\delta}$ of the projection to the second factor.

\begin{claim} \label{section} Let $S=s(\overline{\mathbb D}_{\delta})$ be a  section over
$\overline{\mathbb D}_{\delta}$ such that $S\neq S_u$ and $S\neq S_v$, where $S_u=\{(u,v,\tau)\in \mathbb P^1 \times \overline{\mathbb D}_{\delta}\mid u=0\}$ and $S_v=\{(u,v,\tau)\in \mathbb P^1 \times \overline{\mathbb D}_{\delta}\mid v=0\}$. Then there are functions $f_1(\tau)$ and $f_2(\tau)$ such that
$\text{G.C.D.}(f_1(\tau),f_2(\tau))_{\delta}=1$ and
$f_1(\tau)u+f_2(\tau)v=0$ is an equation of $S$.
\end{claim}
\proof
There are at most finitely many points
$\alpha_1,\dots,\alpha_n$ and $\beta_1,\dots,\beta_m$ in $\overline{\mathbb D}_{\delta}$ such that $s(\alpha_i)\in S_u$, $s(\beta_j)\in S_v$, and $\{ \alpha_1,\dots,\alpha_n\}\cap \{ \beta_1,\dots,\beta_m\}=\emptyset $. Then the section $S$ can be given over $D_u=\overline{\mathbb D}_{\delta}\setminus \{\beta_1,\dots, \beta_m\}$ by equation of the following form $x=h_1(\tau)\frac{(\tau-\alpha_1)^{\mu_1}\dots (\tau-\alpha_n)^{\mu_n} }{(\tau-\beta_1)^{\nu_1}\dots (\tau-\beta_m)^{\nu_m}}$, where $x=\frac{u}{v}$ and $h_1(\tau)$ is a holomorphic invertible function in $D_u$. Similarly, the section $S$ can be given over $D_v=\overline{\mathbb D}_{\delta}\setminus \{\alpha_1,\dots, \alpha_n\}$ by equation of the following form $y=h_2(\tau)\frac{(\tau-\beta_1)^{\nu_1}\dots (\tau-\beta_m)^{\nu_m} }{(\tau-\alpha_1)^{\mu_1}\dots (\tau-\alpha_n)^{\mu_n}}$, where $y=\frac{v}{u}$ and $h_2(\tau)$ is a holomorphic invertible function in $D_v$. In addition, we have $h_1(\tau)=\frac{1}{h_2(\tau)}$ in $D_u\cap D_v$, since $S$ is a  section over
$\overline{\mathbb D}_{\delta}$. Therefore $h_1(\tau)$ can be continued up to the holomorphic function in  $\overline{\mathbb D}_{\delta}$. As a result, we obtain that $S$ can be given by equation ${(\tau-\beta_1)^{\nu_1}\dots (\tau-\beta_m)^{\nu_m}}u-h_1(\tau)(\tau-\alpha_1)^{\mu_1}\dots (\tau-\alpha_n)^{\mu_n} v=0$. \qed \\

Let $h_{\tau}(u,v)$  be a strong equisingular deformation of germs of functions. Denote by $\mathfrak m$  the maximal ideal of the ring $\mathbb C[[u,v]]$ and by $\mu_{\tau_0}(h_{\tau}(u,v))$ the multiplicity of zero of the function $h_{\tau_0}(u,v)$ at the point $(0,0,\tau_0)$, that is, the integer $m$  such that $h_{\tau_0}(u,v)\in \mathfrak m^m\setminus \mathfrak m^{m+1}$.

\begin{claim} \label{proper} Let $h_{\tau}(u,v)$  be a strong equisingular deformation of germs of functions. Then
\begin{itemize}
\item[($i$)] $\mu_{\tau}(h_{\tau}(u,v))$ {\it does not depend on} $\tau$;
\item[($ii$)] {\it if $((h_{\tau_0}(u,v))=B_{\tau_0,1}+\dots +B_{\tau_0,n}$, where $B_{1,\tau_0},\dots,B_{n,\tau_0}$ are the irreducible components of the divisor $((h_{\tau_0}(u,v))$ in $\mathbb D^2_{\varepsilon_1,\varepsilon_2}\times \{ \tau=\tau_0\}$, then $((h_{\tau}(u,v))=\mathcal B_1+\dots +\mathcal B_n$, where $\mathcal B_1,\dots,\mathcal B_n$ are the irreducible components of the divisor $((h_{\tau}(u,v))$ and there are strong equisingular deformations $h_{i,\tau}(u,v)$, $i=1,\dots, n$, such that $((h_{i,\tau}(u,v))=\mathcal B_i$ and $h_{\tau}(u,v)=h_{1,\tau}(u,v)\dots h_{n,\tau}(u,v)$; the functions $h_{i,\tau}(u,v)$ are defined uniquely up to invertible in $\mathcal V$ functions. }
\item[($iii$)] {\it the intersection numbers $\beta_{i,j}(\tau)=(B_{i,\tau},B_{j,\tau})_{o_{\tau}}$, $i\neq j$, at the points $o_{\tau}=(0,0,\tau)$ in $\overline{\mathbb D}^2_{\varepsilon_1,\varepsilon_2}$ do not depend on $\tau\in \overline{\mathbb D}_{\delta}$, where $B_{i,\tau}=\mathcal B_i\cap \text{pr}_2^{-1}(\tau)$.}
\end{itemize}
\end{claim}

\noindent {\it Proof} directly follows from Definition \ref{def1}. \qed\\

Let $u_{1}=u_1(u,v,\tau)$ and $v_{1}=v_1(u,v,\tau)$ be two germs of functions defined in a neighbourhood of the closure of
${\mathbb D}^2_{\varepsilon_1,\varepsilon_2}\times \overline{\mathbb D}_{\delta}\subset \mathbb C^3$ and such that $u_1(0,0,\tau)=0$, $v_1(0,0,\tau)=0$ for all $\tau\in \overline{\mathbb D}_{\delta}$. Assume that $$J(u_{1},v_{1}):= \det \left(\begin{array}{cc}
\frac{\partial u_{1}}{\partial u} & \frac{\partial v_{1}}{\partial u} \\ \frac{\partial u_{1}}{\partial v} & \frac{\partial v_{1}}{\partial v}\end{array}\right)\neq 0 $$ at $(0,0,\tau)$ for all $\tau\in \overline{\mathbb D}_{\delta}$.
Then, by the Inverse Function Theorem,  the triple $(u_1,v_1,\tau)$ is coordinates in some neighbourhood
$V\times\overline{\mathbb D}_{\delta}\subset\mathbb D^2_{\varepsilon_1,\varepsilon_2}\times\overline{\mathbb D}_{\delta}$,
$$V\simeq \mathbb D^2_{\varepsilon'_1,\varepsilon'_2}=\{ (u_1,v_1)\in\mathbb C^2\mid |u_1|<\varepsilon'_1, |v_1|<\varepsilon'_2\},$$ and we will say that the strong equisingular deformation $h_{\tau}(u(u_1,v_1,\tau),v(u_1,v_1,\tau))$ of germs of functions in $\mathbb D^2_{\varepsilon'_1,\varepsilon'_2}$ is obtained from  a strong equisingular deformation  $h_{\tau}(u,v)$ of germs of functions in $\mathbb D^2_{\varepsilon_1,\varepsilon_2}$ by the  coordinate change $(u_1(u,v,\tau), v_1(u,v,\tau))$.

\begin{claim} \label{triv1} Let $h_{\tau}(u,v)$  be a strong equisingular deformation of germs of functions such that $\mu_{\tau}(h_{\tau}(u,v))=1$. Then there is a coordinate change $(u_1(u,v,\tau), v_1(u,v,\tau))$ such that $h_{\tau}(u(u_1,v_1,\tau),v(u_1,v_1,\tau))=u_1$.
\end{claim}
\proof Let $ \displaystyle h_{\tau}(u,v)=\sum_{n=1}^{\infty}\sum_{i+j=n}a_{i,j}(\tau)u^iv^{j} .$
If $a_{1,0}(\tau)\equiv 0$ (the case when $a_{0,1}(\tau)\equiv 0$ is similar) then $a_{0,1}(\tau)$  is an invertible function in $\overline{\mathbb D}_{\delta}$, since $\mu_{\tau}(h_{\tau}(u,v))=1$. Therefore $u_1=h_{\tau}(u,v)$, and $v_1=v$ is a desired coordinate change.

If $a_{0,1}(\tau)a_{1,0}(\tau)\not\equiv 0$ then $\text{G.C.D.}(a_{0,1}(\tau),a_{1,0}(\tau))_{\delta}=1$, since $\mu_{\tau}(h_{\tau}(u,v))=1$.
By Claim \ref{section}, there are functions $b_0(\tau)$ and $b_1(\tau)$ such that $a_{0,1}(\tau)b_1(\tau)-a_{1,0}(\tau)b_0(\tau)=1$. Therefore
$u_1=h_{\tau}(u,v)$,  $v_1=b_1(\tau)u+b_0(\tau)v$  is a desired coordinate change. \qed \\

The proof of the following Proposition essentially repeats the arguments used in \cite{Ar} in the proof of rigidity of curve germs having $ADE$ singularity types, but for completeness it will be given in full.
\begin{prop} \label{ADE1} For any strong equisingular deformation $h_{\tau}(u,v)$ of a function of one of the following singularity types: $A_n$, $n\geq 1$, $D_n$, $n\geq 4$, $E_6$, $E_7$, $E_8$, there is a coordinate change $(\overline u(u,v,\tau), \overline v(u,v,\tau))$ such that the strong equisingular deformation $h_{\tau}(u(\overline u,\overline v,\tau),v(\overline u,\overline v,\tau))$ is trivial.
\end{prop}
\proof Let us write down   the strong equisingular deformation $h_{\tau}(u,v)$ as absolutely converging in ${\mathbb D}^2_{\varepsilon_1,\varepsilon_2}\times \overline{\mathbb D}_{\delta}$ power series
\begin{equation} \label{sed} \displaystyle h_{\tau}(u,v)=\sum_{n=1}^{\infty}\sum_{i+j=n}a_{i,j}(\tau)u^iv^{j} .
\end{equation}

In the beginning, we consider the case where the singularity type is $A_k$, $k\leq 1$. Note that Proposition \ref{ADE1} in the case $k=0$ directly follows from Claim \ref{triv1}.

If $k=1$, then, by Claim \ref{proper}, we have $h_{\tau}(u,v)=h_{1,\tau}(u,v)h_{2,\tau}(u,v)$, where $h_{1,\tau}(u,v)$ and $h_{2,\tau}(u,v)$ are strong equisingular deformations of germs of functions such that $\mu_{\tau}(h_{1,\tau}(u,v))=\mu_{\tau}(h_{2,\tau}(u,v))=1$ and $\beta_{1,2}(\tau)=(B_{1,\tau},B_{2,\tau})_o=1$ for each $\tau\in \overline{\mathbb D}_{\delta}$, where $B_{i,\tau}=(h_{i,\tau}(u,v))$. Therefore, after the coordinate change
\begin{equation} \label{cange1}
\overline u= \frac{1}{2}[h_{1,\tau}(u,v)+h_{2,\tau}(u,v)],\qquad \overline v= \frac{1}{2}[h_{1,\tau}(u,v)-h_{2,\tau}(u,v)],
\end{equation}
we obtain that $h_{\tau}(u(\overline u,\overline v),v(\overline u,\overline v))=\overline u^2-\overline v^2$.

{\bf Case $A_{k}$, $k\geq 2$.} In this case, $\mu_{\tau}(h_{\tau}(u,v))=2$ and, by Proposition \ref{Zar} after $\sigma$-process with center at $o_{\tau}$, the intersection of the proper inverse image of the germ $B$ and the exceptional curve of $\sigma$-process consists of a single point.
Therefore the quadratic form
\begin{equation} \label{Q1} Q_{\tau}(u,v)=a_{2,0}(\tau)u^2+a_{1,1}(\tau)uv+a_{0,2}(\tau)v^2 \end{equation}
is non-trivial at each point $\tau\in\overline{\mathbb D}_{\delta}$ and the discriminant $a_{1,1}^2(\tau)-a_{2,0}(\tau)a_{0,2}(\tau)$ of $Q_{\tau}(u,v)$ is identically equal to zero.
Therefore the equality $Q_{\tau}(u,v)=0$ defines a divisor $2S$ in $\mathbb P^1\times \overline{\mathbb D}_{\delta}$, where $(u,v)$ are homogeneous coordinates in $\mathbb P^1$ and $S=q(\overline{\mathbb D}_{\delta})$ is a section of a holomorphic map
$q:\overline{\mathbb D}_{\delta}\to \mathbb P^1\times\overline{\mathbb D}_{\delta}$ given by solutions
$q(\tau)$  of the equation $Q_{\tau}(u,v)=0$ for each $\tau\in\overline{\mathbb D}_{\delta}$. By Claim \ref{section}, there exist functions $f_1(\tau)$ and $f_2(\tau)$ such that
$\text{G.C.D.}(f_1(\tau),f_2(\tau))_{\delta}=1$ and
$f_1(\tau)u+f_2(\tau)v=0$ is an equation of $S$. Therefore $Q_{\tau}(u,v)=f_0(\tau)(f_1(\tau)u+f_2(\tau)v)^2$, where $f_0(\tau)$ is an invertible in $\overline{\mathbb D}_{\delta}$ function.

Applying Claim \ref{triv1} to the strong equisingular deformation $$\sqrt{Q_{\tau}(u,v)}=\sqrt{f_0(\tau)}(f_1(\tau)u+f_2(\tau)v),$$ we find a coordinate change
$u(u_0,v_0,\tau)$, $v(u_0,v_0,\tau)$ such that  $h_{\tau}(u(u_0,v_0),v(u_0,v_0))$ has the following form:
$$h_{\tau}(u(u_0,v_0),v(u_0,v_0))=u_0^2+\sum_{m=3}^{\infty}\,\, \sum_{i+j=m}a_{0,i,j}(\tau)u_0^iv_0^j.$$

After the coordinate change:
$$u_1=u_0\sqrt{1+\sum_{m=3}^{\infty}a_{0,m,0}(\tau)u_0^{m-2}},\quad v_1=v_0$$
we obtain that $h_{\tau}(u(u_1,v_1),v(u_1,v_1))$ has the following form:
$$h_{\tau}(u(u_1,v_1),v(u_1,v_1))=u_1^2-v_1\sum_{m=m_0}^{\infty}\,\, \sum_{i=0}^ma_{1,i,m-i+1}(\tau)u_1^iv_1^{m-i}$$
(here $m_0=2$, but below we will assume that $m_0=1$ to include in consideration the case of singularities of type $A_1$).

Next, we consistently make $k$  changes of coordinates of the following form:
\begin{equation} \label{seq1} u_{l+1}=u_l-\frac{1}{2}a_{l,1,l}(\tau)v_l^l,\quad v_{l+1}=v_l.
\quad \displaystyle l=1,\dots,k.
\end{equation}
Since the functions $h_{\tau}(u(u_{l+1},v_{l+1}),v(u_{l+1},v_{l+1}))$ are strong equisingular deformation of germs of functions of singularity type $A_k$, it follows from Proposition \ref{Zar} that
\begin{equation} \label{seq2} \begin{array}{c} h_{\tau}(u(u_{l+1},v_{l+1}),v(u_{l+1},v_{l+1}))=  u_{l+1}^2[1- \displaystyle
\sum_{m=3}^{\infty}\,\,\sum_{i=2}^ma_{l+1,i,m-i}(\tau)u_{l+1}^{i-2}v_{l+1}^{m-i}]-  \\
\displaystyle -\sum_{m=l+1}^{\infty}[a_{l+1,0,m}(\tau)+a_{l+1,1,m}(\tau)u_{l+1}]v_{l+1}^{m},
\end{array}\end{equation}
where if $l=k$ then $a_{k+1,0,k+1}(\tau)$ is an invertible function.
 Therefore, after the coordinate change:
$$\begin{array}{ll} \overline u= & u_{k+1}\sqrt{1- \displaystyle
\sum_{m=3}^{\infty}\,\,\sum_{i=2}^ma_{2,i,m-i}(\tau)u_{k+1}^{i-2}v_{k+1}^{m-i}},\\ \overline v= & v_{k+1}\sqrt[k+1]{\displaystyle \sum_{m=k+1}^{\infty}[a_{2,0,m}(\tau)+a_{2,1,m}(\tau)u_{k+1}]v_{k+1}^{m-k-1}},\end{array}$$
we obtain that  $h_{\tau}(u(\overline u,\overline v),v(\overline u,\overline v))=\overline u^2-\overline v^{k+1}$.

{\bf Case $D_n$, $n\geq 4$.} Put $k=n-3$. By Claim \ref{proper}, we have
$h_{\tau}(u,v)=h_{1,\tau}(u,v)h_{2,\tau}(u,v)$, where $h_{1,\tau}(u,v)$ is a strong equisingular deformation of germs of functions such that $\mu_{\tau}(h_{1,\tau}(u,v))=1$, and $h_{2,\tau}(u,v)$ is a  strong equisingular deformation of germs of functions of singularity type $A_{k}$, and $\beta_{1,2}(\tau)=(B_{1,\tau},B_{2,\tau})_{o_{\tau}}=2$ for each $\tau\in \overline{\mathbb D}_{\delta}$, where $B_{i,\tau}=(h_{i,\tau}(u,v))$. Therefore there is a coordinate change $(u_0(u,v,\tau), v_0(u,v,\tau))$ such that
$h_{1,\tau}(u(u_0,v_0),v(u_0,v_0))=v_0$ and
$$h_{2,\tau}(u(u_0,v_0),v(u_0,v_0))=\sum_{m=2}^{\infty}\,\,\sum_{i+j=m}a_{0,i,j}(\tau)u_0^iv_0^j,$$
where $a_{2,0}(\tau)$ is an invertible in $\overline{\mathbb D}_{\delta}$ function, since  $\beta_{1,2}(\tau)
=2$ for each $\tau\in \overline{\mathbb D}_{\delta}$.  It is directly checked that after the change of coordinates
$$u_1=u_0\sqrt{\sum_{i=2}^{\infty}a_{0,i,0}(\tau)u_0^{i-2}}, \quad v_1=v_0, $$
the function $h_{\tau}(u(u_1,v_1),v(u_1,v_1))$ has the following form:
$$h_{\tau}(u(u_1,v_1),v(u_1,v_1))=v_1[u_1^2-
v_1\sum_{m=1}^{\infty}\,\sum_{i+j=m}a_{0,i,j}(\tau)u_1^iv_1^j].$$
Next, as in the case of $A_k$, we consistently make $k$ coordinate changes  $u=u(u_l,v_l)), v=v(u_l,v_l)$, $l=1,\dots, k$ (see (\ref{seq1})),
and write down the deformation $h_{\tau}(u,v)$ as a function of variables $u_l$ and $v_l$ (see (\ref{seq2})). The result is that
$$\begin{array}{r} h_{\tau}(u(u_{k+1},v_{k+1}),v(u_{k+1},v_{k+1}))=  v_{k+1}[ u_{k+1}^2(1- \displaystyle
\sum_{m=3}^{\infty}\,\,\sum_{i=2}^ma_{k+1,i,m-i}(\tau)u_{k+1}^{i-2}v_{k+1}^{m-i})-  \\
-\displaystyle v_{k+1}^{k+1}[\sum_{m=k+1}^{\infty}(a_{k+1,0,m}(\tau)+
 a_{k+1,1,m}(\tau))u_{k+1}v_{k+1}^{m-k-1}] ,\end{array}$$
where  $a_{k+1,0,k+1}(\tau)$ is an invertible function.

Finally, after the coordinate change
$$\overline u=f(u_{k+1},v_{k+1},\tau)u_{k+1},\quad \overline v=g(u_{k+1},v_{k+1},\tau)v_{k+1},$$
where
$$\begin{array}{ll} \\     g(u_{k+1},v_{k+1},\tau)= &
\displaystyle \sqrt[k+2]{\sum_{j=k+1}^{\infty}(a_{k+1,1,j}(\tau)u_{k+1}+  a_{k+1,0,j}(\tau))v_{k+1}^{j-k-1}}, \\ f(u_{k+1},v_{k+1},\tau)= &
\displaystyle \sqrt{g(u_{k+1},v_{k+1},\tau)^{-1}[1+ \displaystyle \sum_{m=0}^{\infty}\,\,\sum_{i+j=m}a_{k+1,i,j}(\tau)u_{k+1}^iv_{k+1}^{j}]}, \end{array}$$
we obtain that $h_{\tau}(u(\overline u,\overline v),v(\overline u,\overline v))=  \overline v[\overline u^2-\overline v^{k+1}]$.

{\bf Case $E_7$.} By Claim \ref{proper}, we have
$h_{\tau}(u,v)=h_{1,\tau}(u,v)h_{2,\tau}(u,v)$, where $h_{1,\tau}(u,v)$ is a strong equisingular deformation of germs of functions such that $\mu_{\tau}(h_{1,\tau}(u,v))=1$, and $h_{2,\tau}(u,v)$ is a  strong equisingular deformation of germs of functions of singularity type $A_{2}$, and $\beta_{1,2}(\tau)=(B_{1,\tau},B_{2,\tau})_o=3$ for each $\tau\in \overline{\mathbb D}_{\delta}$, where $B_{i,\tau}=(h_{i,\tau}(u,v))$. Therefore there is a coordinate change $(u_0(u,v,\tau), v_0(u,v,\tau))$ such that
$h_{1,\tau}(u(u_0,v_0),v(u_0,v_0))=u_0$ and
$$h_{2,\tau}(u(u_0,v_0),v(u_0,v_0))=\sum_{m=2}^{\infty}\,\,\sum_{i+j=m}a_{0,i,j}(\tau)u_0^iv_0^j,$$
where $a_{0,2,0}(\tau)$ and $a_{0,0,3}(\tau)$
are invertible in $\overline{\mathbb D}_{\delta}$ functions and $a_{0,1,1}(\tau)=a_{0,0,2}(\tau)\equiv 0$, since the quadratic homogeneous  form
$a_{0,2,0}(\tau)u_0^2+a_{0,1,1}(\tau)u_0v_0+a_{0,0,2}(\tau)v_0^2$ is the square of a linear non-generate form for a singularity of type $A_2$ and $\beta_{1,2}(\tau)
=3$ for each $\tau\in \overline{\mathbb D}_{\delta}$.  It is directly checked that after the change of coordinates
$$u_1=\sqrt[3]{a_{0,2,0}(\tau)}u_0, \quad v_1=v_0$$
we obtain
$$h_{\tau}(u(u_1,v_1),v(u_1,v_1))=u_1^3-
u_1\sum_{m=3}^{\infty}\,\,\sum_{i=0}^ma_{1,i,m-i}(\tau)u_1^iv_1^{m-i},$$
where $a_{1,0,3}(\tau)$ is an invertible function. After the coordinate change
$$u_2=u_1-\frac{1}{3}a_{1,1,2}v_1^2, \quad v_2=v_1,$$
we obtain
$$\begin{array}{ll} h_{\tau}(u(u_2,v_2),v(u_2,v_2))= & \displaystyle u_2^3[1-\sum_{i=3}^{\infty}
(a_{2,i,0}u_2^{i-2}+a_{2,i,1}(\tau)u_2^{i-2}v_2+a_{2,i,2}(\tau)u_2^{i-2}v_2^2]- \\ &
\displaystyle -u_2v_2^3[a_{2,0,3}(\tau)+\sum_{m=4}^{\infty}\,\,\sum_{i+j=m}a_{2,i,j}(\tau)u_2^iv_2^{j-3}]\end{array},$$
where $a_{2,0,3}(\tau)$ is an invertible function in $\overline{\mathbb D}_{\delta}$. Finally, after the coordinate change
$$ \begin{array}{ll} \overline u= & \displaystyle u_2\sqrt[3]{1-\sum_{i=3}^{\infty}
[a_{2,i,0}u_2^{i-2}+a_{2,i,1}(\tau)u_2^{i-2}v_2+a_{2,i,2}(\tau)u_2^{i-2}v_2^2]}, \\
\overline v= & v_2\frac{\sqrt[3]{\displaystyle a_{2,0,3}(\tau)+  \sum_{m=4}^{\infty}\,\,  \sum_{i+j=m}a_{2,i,j}(\tau)u_2^iv_2^{j-3}}}{\sqrt[9]{ \displaystyle 1-  \sum_{i=3}^{\infty}
[a_{2,i,0}u_2^{i-2}+a_{2,i,1}(\tau)u_2^{i-2}v_2+a_{2,i,2}(\tau)u_2^{i-2}v_2^2]}} \end{array} $$
we obtain that $h_{\tau}(u(\overline u,\overline v),v(\overline u,\overline v))=\overline u(\overline u^2-\overline v^3)$.

{\bf Case $E_6$.} In this case $\mu_{\tau_0}(h_{\tau_0}(u,v))=3$  for all $\tau_0\in \overline{\mathbb D}_{\delta}$ and by Proposition \ref{Zar} after the $\sigma$-process with center at $o_{\tau_0}$, the proper inverse image $B'_{\tau_0}$ of the irreducible germ $B_{\tau_0}=(h_{\tau_0}(u,v))$ is non-singular and the intersection number $(E,B'_{\tau_0})$ of the germ $B'_{\tau_0}$ and the exceptional curve $E$ 
is equal to three.
Therefore the cubic homogeneous form
$$ C_{\tau}(u,v)=\sum_{i+j=3}a_{i,j}(\tau)u^iv^j$$
is the cube of a  linear form $L_{\tau}(u,v)=b_{1,0}(\tau)u+b_{0,1}(\tau)v$ non-degenerated at each point $\tau\in \overline{\mathbb D}_{\delta}$. By Claim \ref{line}, there is a coordinate change $u_0=L_{\tau}(u,v)$ and $v_0=c_{1,0}(\tau)u+c_{0,1}(\tau)v$ such that $h_{\tau}(u(u_0,v_0),v(u_0,v_0))$ has the following form:
$$h_{\tau}(u(u_0,v_0),v(u_0,v_0))= u_0^3- \sum_{m=4}^{\infty}\,\,\sum_{i+j=m}a_{0,i,j}(\tau)u_0^iv_0^j.$$

Next, we make the change of coordinates of the following form:
$$\displaystyle u_1=u_0-\frac{1}{3}\sum_{i=2}^4a_{0,i,4-i}u_0^{i-2}v_0^{4-i},\quad v_1=v_0$$ 
and we obtain that 
\begin{equation} \label{seq4}  h_{\tau}(u(u_{1},v_{1}),v(u_{1},v_{1}))=  u_{1}^3-a_{1,1,3}(\tau)u_1v_1^3- a_{1,0,4}(\tau)v_1^4-
\displaystyle  \sum_{m=5}^{\infty}\,\,\sum_{i+j=m}a_{1,i,j}(\tau)u_{1}^iv_{1}^{j},\end{equation}
where the function $a_{1,0,4}(\tau)$ is invertible, since $B'_{\tau}$ is an irreducible germ of a non-singular curve and $(E,B'_{\tau})=3$.
After the coordinate change
$$u_{2}= u_1, \quad v_2= \sqrt[4]{a_{1,0,4}(\tau)}v_1 +\frac{a_{1,1,3}}{4\sqrt[4]{a_{1,0,4}^3(\tau)}} u_{1}, $$
we obtain that
$$h_{\tau}(u(u_{2},v_{2}),v(u_{2},v_{2}))= \displaystyle u_{1}^3-u_2^2\sum_{i=2}^4a_{2,i,4-i}(\tau)u_2^{i-2}v_2^{4-i}-
v_2^4-\displaystyle  \sum_{m=5}^{\infty}\,\,\sum_{i+j=m}a_{2,i,j}(\tau)u_{2}^iv_{2}^{j}. $$
After that we make
the change of coordinates
$$\displaystyle u_3=u_2-\frac{1}{3}\sum_{i=2}^4a_{2,i,4-i}u_2^{i-2}v_2^{4-i},\quad v_3=v_2  $$
and we obtain that $h_{\tau}(u(u_{3},v_{3}),v(u_{3},v_{3}))$ has the following form:
$$h_{\tau}(u(u_{3},v_{3}),v(u_{3},v_{3}))=  u_{3}^3-v_3^4 -\displaystyle  \sum_{m=5}^{\infty}\,\,\sum_{i+j=m}a_{3,i,j}(\tau)u_{3}^iv_{3}^{j}. $$
After
the change of coordinates
$\displaystyle u_4=u_3-\frac{1}{3}a_{3,2,3}v_2^{3},\,\, v_4=v_3,  $
we obtain that
$$\begin{array}{ll} h_{\tau}(u(u_{4},v_{4}),v(u_{4},v_{4}))= &
\displaystyle [u_4^3-\sum_{m=5}^{\infty}\,\,\sum_{i=3}^ma_{4,i,m-i}(\tau)u_4^{i}v_4^{m-i}]-  \\
& \displaystyle 
-[v_4^4+\sum_{m=5}^{\infty}\,\,\sum_{i=0}^2a_{4,i,m-i}(\tau)u_{4}^iv_{4}^{m-i}], \end{array}
$$
where $a_{4,2,3}(\tau)\equiv 0$. Finally, after the coordinate change
$$\begin{array}{ll} \overline u= & u_4\sqrt[3]{1-\sum_{m=5}^{\infty}\,\,\sum_{i=3}^ma_{4,i,m}(\tau)u_4^{i-3}v_4^{m-i}},\\
\overline v= & v_4\sqrt[4]{1+\sum_{m=5}^{\infty}\,\,\sum_{i=0}^2a_{4,i,m-i-4}(\tau)u_{4}^iv_{4}^{m-i-4}}, \end{array}$$
we obtain that $h_{\tau}(u(\overline u,\overline v),v(\overline u,\overline v))=\overline u^3-\overline v^4$.

{\bf Case $E_8$.}  After   the coordinate change $(u(u_1,v_1,\tau),v(u_1,v_1,\tau)$, the same as in the case $E_6$, we obtain that
$h_{\tau}(u(u_1,v_1),v(u_1,v_1)$ can be written in form (\ref{seq4}) in which $a_{1,0,4}(\tau)=a_{1,1,3}(\tau)\equiv 0$, since
the proper inverse image $B_{1,\tau_0}$ of the divisor $B_{\tau_0}=(h_{\tau_0}(u,v))$ after the 
$\sigma$-process at $o_{\tau_0}$ has the singularity type $A_2$ for each
$\tau_0\in \overline{\mathbb D}_{\delta}$.

After the coordinate change $u_2=u_1-\frac{1}{3}\sum_{i=2}^5a_{1,i,5-i}u_1^{i-2}v_1^{5-i}$, $v_2=v_1$, we obtain
$$h_{\tau}(u(u_2,v_2),v(u_2,v_2))=u_2^3-a_{2,0,5}(\tau)v_2^5-a_{2,1,4}(\tau)u_2v_2^4- \displaystyle \sum_{m=6}^{\infty}\,\, \sum_{i=0}^ma_{2,i,m-i}(\tau)u_2^iv_2^{m-i},$$
where $a_{2,0,5}(\tau)$ is an invertible function, since the divisors $B_{1,\tau_0}$ have the singularity of type $E_8$ for all $\tau_0\in \overline{\mathbb D}_{\delta}$. Therefore, after the coordinate change
$$ u_3=u_2,\quad v_3=\sqrt[5]{a_{2,0,5}(\tau)}v_2+\frac{a_{2,1,4}(\tau)}{5\sqrt[5]{a_{2,0,5}(\tau)^4}}u_2,$$
we obtain
$$h_{\tau}(u(u_3,v_3),v(u_3,v_3))=u_3^3-v_3^5 \displaystyle -\sum_{i=2}^5a_{3,i,5-i}(\tau)u_3^iv_3^{5-i}-  \sum_{m=6}^{\infty}\,\, \sum_{i=0}^ma_{3,i,m-i}(\tau)u_3^iv_3^{m-i}.$$
Now, after the coordinate change
$$u_4=u_3 \displaystyle -\frac{1}{3}\sum_{i=2}^5a_{3,i,5-i}(\tau)u_3^{i-2}v_3^{5-i}, \quad v_4=v_3,$$
we obtain
$$h_{\tau}(u(u_4,v_4),v(u_4,v_4))=u_4^3-v_4^5 \displaystyle -  \sum_{m=6}^{\infty}\,\, \sum_{i=0}^ma_{4,i,m-i}(\tau)u_4^iv_4^{m-i}.$$
The coordinate change $u_5=\frac{1}{3}a_{4,2,4}(\tau)v_4^4$, $v_5=v_4$ gives
$$\begin{array}{ll} h_{\tau}(u(u_5,v_5),v(u_5,v_5))= & u_5^3[1 \displaystyle -  \sum_{m=6}^{\infty}\,\, \sum_{i=3}^ma_{4,i,m-i}(\tau)u_5^{i-3}iv_4^{m-i}]- \\ \phantom{ h_{\tau}(u(u_5,v_5),v(u_5,v_5))} &
-v_5^5[1 \displaystyle +  \sum_{m=6}^{\infty}\,\, \sum_{i=0}^2a_{5,i,m-i-5}(\tau)u_5^iv_5^{m-i-5}],\end{array}$$
where $a_{5,2,-1}(\tau)\equiv 0$. Finally, after the coordinate change
$$ \begin{array}{ll}
\overline u= & u_5\sqrt[3]{1 \displaystyle -  \sum_{m=6}^{\infty}\,\, \sum_{i=3}^ma_{4,i,m-i}(\tau)u_5^{i-3}iv_4^{m-i}}, \\
\overline v = &
v_5\sqrt[5]{1 \displaystyle +  \sum_{m=6}^{\infty}\,\, \sum_{i=0}^2a_{5,i,m-i-5}(\tau)u_5^iv_5^{m-i-5}}, \end{array}
$$
we obtain that $h_{\tau}(u(\overline u,\overline v),v(\overline u,\overline v))=\overline u^3-\overline v^5$. \qed \\

Obviously, to prove Theorem \ref{main-t} it suffices to  prove
\begin{prop} \label{Prop2} Let $\mathcal F=F_{\tau}:(U,o')\times \overline{\mathbb D}_{\delta}\to
\mathbb D^2_{\varepsilon_1,\varepsilon_2}\times\overline{\mathbb D}_{\delta}$ be a strong deformation of a germ of cover
$F_0: (U,o')\times\{ \tau=0\}\to \mathbb D^2_{\varepsilon_1,\varepsilon_2}\times\{ \tau=0\}$ branched along a divisor $B_0$ having one of $ADE$ singularity types. Then for each $\tau_0\in\overline{\mathbb D}_{\delta}$, the germs of covers $F_{\tau_0}: (U,o')\times\{ \tau=\tau_0\}\to \mathbb D^2_{\varepsilon_1,\varepsilon_2}\times\{ \tau=\tau_0\}$
are equivalent to the germ $F_0$.
\end{prop}
\proof Let $h_{\tau}(u,v)=0$ be an equation of the branch divisor $\mathcal B$ of the strong deformation $\mathcal F$.
By Proposition \ref{ADE1}, there is a coordinate change $(\overline u(u,v,\tau),\overline v(u,v,\tau),\tau)$ in a neighbourhood
$ \mathcal V 
\subset \mathbb D^2_{\varepsilon_1,\varepsilon_2}\times\overline{\mathbb D}_{\delta}$, where $ \mathcal V \simeq V\times
\overline{\mathbb D}_{\delta}$ and
$$V\simeq \mathbb D^2_{\varepsilon'_1,\varepsilon'_2}=\{ (\overline u,\overline v)\in \mathbb C^2\mid |\overline u|<\varepsilon'_1,|\overline v|<\varepsilon'_2\},$$
in which the function $h_{\tau}(u(\overline u,\overline v),v(\overline u,\overline v))$ does not depend on $\tau$.

The coordinate change $(u,v,\tau)\mapsto(\overline u(u,v,\tau),\overline v(u,v,\tau),\tau)$ defines a bi-holomorphic mapping
$$\mathcal G:\mathcal V\to \mathbb D^2_{\varepsilon'_1,\varepsilon'_2}\times\overline{\mathbb D}_{\delta}, \quad
(u,v,\tau)\mapsto (\overline u(u,v,\tau),\overline v(u,v,\tau),\tau).$$
Denote by
$$\mathcal W=\mathcal F^{-1}(\mathcal V),\,\,\,\,\, \mathcal H:= \mathcal G\circ \mathcal F:\mathcal W \to \mathbb D^2_{\varepsilon'_1,\varepsilon'_2}\times\overline{\mathbb D}_{\delta},$$ 
$$V_{\tau_0}=V\times \{\tau=\tau_0\},\,\,\,\,\,    W_{\tau_0}=\mathcal F^{-1}(\mathcal G^{-1}(V_{\tau_0}\})).$$

Obviously, for each $\tau_0\in \overline{\mathbb D}_{\delta}$ the covers $H_{\tau_0}:W_{\tau_0}\to V_{\tau_0}$ and  $F_{\tau_0}:W_{\tau_0}\to \mathcal G^{-1}(V_{\tau_0})$ are equivalent.

The cover $\mathcal H:\mathcal W \to \mathbb D^2_{\varepsilon'_1,\varepsilon'_2}\times\overline{\mathbb D}_{\delta}$ is branched in $\mathcal G(\mathcal B)=B\times \overline{\mathbb D}_{\delta}$, where $B=\mathcal G(B_0\cap\mathcal V)$, and it induces a monodromy homomorphism
$\mathcal H_*:\pi_1((\mathbb D^2_{\varepsilon'_1,\varepsilon'_2}\times\overline{\mathbb D}_{\delta})\setminus (B\times \overline{\mathbb D}_{\delta}))\to \mathbb S_d$.

Similarly, for each $\tau_0\in \overline{\mathbb D}_{\delta}$
the cover $H_{\tau_0}:W_{\tau_0} \to V_{\tau_0}$ is branched in $B\times \{ \tau=\tau_0\}$ and it defines a monodromy homomorphism
$$H_{\tau_0*}=\mathcal H_*\circ i_{\tau_{0*}}:\pi_1(V_{\tau_0}\setminus (B\times \{\tau=\tau_0\})\to \mathbb S_d,$$ where
$i_{\tau_{0_*}}:\pi_1(V_{\tau_0}\setminus (B\times \{\tau=\tau_0\})\to \pi_1((\mathbb D^2_{\varepsilon'_1,\varepsilon'_2}\times\overline{\mathbb D}_{\delta})\setminus (B\times \overline{\mathbb D}_{\delta}))$ is an isomorphism induced by imbedding $i_{\tau_0}:V_{\tau_0}\setminus (B\times \{\tau=\tau_0\})\hookrightarrow (\mathbb D^2_{\varepsilon'_1,\varepsilon'_2}\times\overline{\mathbb D}_{\delta})\setminus (B\times \overline{\mathbb D}_{\delta})$. The identification (due to $\text{pr}_1$) of pairs $(V,B)$ and $(V_{\tau_0},B\times \{ \tau=\tau_0\})$ gives rise to identification of homomorphisms $H_{0*}$ and $H_{\tau_0*}$. Therefore we can identify the covers
$H_0:W_0\setminus H^{-1}_0\to V_0\setminus B$ and $H_{\tau_0}:W_{\tau_0}\setminus H^{-1}_{\tau_0}\to V_0\setminus (B\times\{\tau=\tau_0\})$ and hence, by Grauert - Remmert - Riemann - Stein Theorem (\cite{G-R}), the covers $H_0$ and $H_{\tau_0}$ are equivalent. \qed

\section{Local fundamental groups of curve germs}
Let  $\Gamma(B)$ be the graph of a curve germ $(B,o)$ (below we use notations used in Definition  \ref{grph}). Note that $\Gamma(B)$ is a tree. Let us call the vertex $e_{1+m}$ the {\it root} of $\Gamma(B)$ and renumber the vertices $e{2+m},\dots,e_{n+m}$ (and the corresponding to them curves $E_{i+m}$) so that the new numbering has the  following property:
\newline
$\phantom{aaaa}$ {\it in the shortest path
$(e_{1+m},e_{i_1}), (e_{i_1},e_{i_2}),\dots, (e_{i_{k-1}},e_{i_{k}})$
 from $e_{1+m}$ to each vertex \newline $\phantom{aaaa}$ $e_{i_k}$  along edges $(e_{i_{j-1}},e_{i_{j}})$, 
 we have inequalities
$i_{j-1}<i_{j}$ for $j=1,\dots, k$.}

For  vertices $v_i$ and $v_j$ of the graph $\Gamma(B)$ we define
$$ \delta_{i,j}=\left\{ \begin{array}{ll} 1, & \text{if}\,\, v_i\,\, \text{and}\,\, v_j\,\, \text{are connected by an edge in}\,\, \Gamma(B),
\\ 0, & \text{if}\,\, v_i\,\, \text{and}\,\, v_j\,\, \text{are not connected by an edge in}\,\, \Gamma(B), \\ 0, & \text{if}\,\, i=j.\end{array} \right. $$

The graph $\Gamma(B)$ of the curve germ $B$ of singularity type $A_{2n+1}$, $n\geq 0$, is depicted on Fig. 1 (if $n=0$ then the weight  of the vertex $e_3$ is equal to $-1$).

\begin{picture}(300,85)
\put(85,30){\circle*{3}}\put(80,35){$\mbox{}_{-2}$}
\put(82,20){$\mbox{e}_{3}$}\put(85,30){\line(1,0){40}}
\put(125,30){\circle*{3}}\put(120,35){$\mbox{}_{-2}$}
\put(122,20){$\mbox{e}_{4}$}\put(125,30){\line(1,0){40}}
\put(175,30){$\dots$} \put(275,30){\circle*{2}}
\put(227,35){$\mbox{}_{-2}$}\put(235,30){\circle*{3}}
\put(233,20){$\mbox{e}_{n+2}$}\put(195,30){\line(1,0){40}}
\put(263,35){$\mbox{}_{-1}$}\put(275,30){\circle*{3}}
\put(270,20){$\mbox{e}_{n+3}$}\put(235,30){\line(1,0){40}}
\put(275,70){\circle*{3}}
\put(280,68){$\mbox{b}_{1}$}\put(275,30){\line(0,1){40}}
\put(315,30){\circle*{3}}
\put(313,20){$\mbox{b}_{2}$}\put(275,30){\line(1,0){40}}
\put(200,-8){$\text{Fig.}\, 1$}
\end{picture} \vspace{0.5cm}

The graph $\Gamma(B)$ of the curve germ $B$ of singularity type $A_{2n}$, $n\geq 1$, is depicted on Fig. 2.

\begin{picture}(300,85)
\put(85,30){\circle*{3}}\put(80,35){$\mbox{}_{-2}$}
\put(82,20){$\mbox{e}_{2}$}\put(85,30){\line(1,0){40}}
\put(125,30){\circle*{3}}\put(120,35){$\mbox{}_{-2}$}
\put(122,20){$\mbox{e}_{3}$}\put(125,30){\line(1,0){40}}
\put(175,30){$\dots$} \put(275,30){\circle*{2}}
\put(227,35){$\mbox{}_{-3}$}\put(235,30){\circle*{3}}
\put(233,20){$\mbox{e}_{n+1}$}\put(195,30){\line(1,0){40}}
\put(263,35){$\mbox{}_{-1}$}\put(275,30){\circle*{3}}
\put(270,20){$\mbox{e}_{n+2}$}\put(235,30){\line(1,0){40}}
\put(275,70){\circle*{3}}
\put(280,68){$\mbox{b}_{1}$}\put(275,30){\line(0,1){40}}
\put(307,35){$\mbox{}_{-2}$}\put(315,30){\circle*{3}}
\put(313,20){$\mbox{e}_{n+3}$}\put(275,30){\line(1,0){40}}
\put(200,-8){$\text{Fig.}\, 2$}
\end{picture} \vspace{0.5cm}

The graph $\Gamma(B)$ of the curve germ $B$ of singularity type $D_{2n+2}$, $n\geq 1$, is depicted on Fig. 3.

\begin{picture}(300,85)
\put(85,30){\circle*{3}}
\put(82,20){$\mbox{b}_{1}$}\put(85,30){\line(1,0){40}}
\put(125,30){\circle*{3}}\put(120,35){$\mbox{}_{-2}$}
\put(122,20){$\mbox{e}_{4}$}\put(125,30){\line(1,0){40}}
\put(175,30){$\dots$} \put(275,30){\circle*{2}}
\put(227,35){$\mbox{}_{-2}$}\put(235,30){\circle*{3}}
\put(233,20){$\mbox{e}_{n+2}$}\put(195,30){\line(1,0){40}}
\put(263,35){$\mbox{}_{-1}$}\put(275,30){\circle*{3}}
\put(270,20){$\mbox{e}_{n+3}$}\put(235,30){\line(1,0){40}}
\put(275,70){\circle*{3}}
\put(280,68){$\mbox{b}_{2}$}\put(275,30){\line(0,1){40}}
\put(315,30){\circle*{3}}
\put(313,20){$\mbox{b}_{3}$}\put(275,30){\line(1,0){40}}
\put(200,-8){$\text{Fig.}\, 3$}
\end{picture} \vspace{0.5cm}

The graph $\Gamma(B)$ of the curve germ $B$ of singularity type $D_{2n+3}$, $n\geq 1$, is depicted on Fig. 4.

\begin{picture}(300,85)
\put(85,30){\circle*{3}}
\put(82,20){$\mbox{b}_{1}$}\put(85,30){\line(1,0){40}}
\put(125,30){\circle*{3}}\put(120,35){$\mbox{}_{-2}$}
\put(122,20){$\mbox{e}_{3}$}\put(125,30){\line(1,0){15}}
\put(150,30){$\dots$} \put(275,30){\circle*{2}}
\put(175,30){\line(1,0){15}}\put(193,30){\circle*{3}}\put(192,35){$\mbox{}_{-2}$}
\put(190,20){$\mbox{e}_{n+1}$}
\put(227,35){$\mbox{}_{-3}$}\put(235,30){\circle*{3}}
\put(233,20){$\mbox{e}_{n+2}$}\put(195,30){\line(1,0){40}}
\put(263,35){$\mbox{}_{-1}$}\put(275,30){\circle*{3}}
\put(270,20){$\mbox{e}_{n+3}$}\put(235,30){\line(1,0){40}}
\put(275,70){\circle*{3}}
\put(280,68){$\mbox{b}_{2}$}\put(275,30){\line(0,1){40}}
\put(307,35){$\mbox{}_{-2}$}\put(315,30){\circle*{3}}
\put(313,20){$\mbox{e}_{n+4}$}\put(275,30){\line(1,0){40}}
\put(200,-8){$\text{Fig.}\, 4$}
\end{picture} \vspace{0.5cm}

The graph $\Gamma(B)$ of the curve germ $B$ of singularity type $E_{6}$ is depicted on Fig. 5.

\begin{picture}(300,85)
\put(147,35){$\mbox{}_{-4}$}\put(155,30){\circle*{3}}
\put(153,20){$\mbox{e}_{2}$}
\put(183,35){$\mbox{}_{-1}$}\put(195,30){\circle*{3}}
\put(190,20){$\mbox{e}_{3}$}\put(155,30){\line(1,0){40}}
\put(195,70){\circle*{3}}
\put(200,68){$\mbox{b}_{1}$}\put(195,30){\line(0,1){40}}
\put(227,35){$\mbox{}_{-2}$}\put(235,30){\circle*{3}}
\put(233,20){$\mbox{e}_{4}$}\put(195,30){\line(1,0){40}}
\put(267,35){$\mbox{}_{-2}$}\put(275,30){\circle*{3}}
\put(273,20){$\mbox{e}_{5}$}\put(235,30){\line(1,0){40}}
\put(200,-8){$\text{Fig.}\, 5$}
\end{picture} \vspace{0.5cm}

The graph $\Gamma(B)$ of the curve germ $B$ of singularity type $E_{7}$ is depicted on Fig. 6.

\begin{picture}(300,85)
\put(127,35){$\mbox{}_{-3}$}\put(135,30){\circle*{3}}
\put(132,20){$\mbox{e}_{3}$}
\put(163,35){$\mbox{}_{-1}$}\put(175,30){\circle*{3}}
\put(170,20){$\mbox{e}_{4}$}\put(135,30){\line(1,0){40}}
\put(175,70){\circle*{3}}
\put(180,68){$\mbox{b}_{2}$}\put(175,30){\line(0,1){40}}
\put(207,35){$\mbox{}_{-2}$}\put(215,30){\circle*{3}}
\put(213,20){$\mbox{e}_{5}$}\put(175,30){\line(1,0){40}}
\put(255,30){\circle*{3}}
\put(253,20){$\mbox{b}_{1}$}\put(215,30){\line(1,0){40}}
\put(200,-8){$\text{Fig.}\, 6$}
\end{picture} \vspace{0.5cm}

The graph $\Gamma(B)$ of the curve germ $B$ of singularity type $E_{8}$ is depicted on Fig. 7.

\begin{picture}(300,85)
\put(127,35){$\mbox{}_{-3}$}\put(135,30){\circle*{3}}
\put(132,20){$\mbox{e}_{5}$}
\put(163,35){$\mbox{}_{-1}$}\put(175,30){\circle*{3}}
\put(170,20){$\mbox{e}_{4}$}\put(135,30){\line(1,0){40}}
\put(175,70){\circle*{3}}
\put(180,68){$\mbox{b}_{1}$}\put(175,30){\line(0,1){40}}
\put(207,35){$\mbox{}_{-2}$}\put(215,30){\circle*{3}}
\put(213,20){$\mbox{e}_{3}$}\put(175,30){\line(1,0){40}}
\put(247,35){$\mbox{}_{-3}$}\put(255,30){\circle*{3}}
\put(253,20){$\mbox{e}_{2}$}\put(215,30){\line(1,0){40}}
\put(200,-8){$\text{Fig.}\, 7$}
\end{picture} \vspace{0.5cm}

\begin{rem} \label{rem1} Note that in all graphs $\Gamma(B)$ of curve germs $B$ of $ADE$ singularity types (except the singularity types $A_0$ and $A_1$) there is a single vertex of valency three and this vertex has weight $w=-1$. \end{rem}

For a curve germ $B$, denote by the same letters $b_i$ elements in the local fundamental group $\pi_1^{loc}(B,o)=\pi _1(V\setminus B, p)\simeq \pi _1(V_n\setminus \sigma^{-1}(B), \sigma^{-1}(p))$  represented by some loops  $\lambda_j\subset V_n\setminus \sigma^{-1}(B)$ around $B'_j$ and by $e_i$ elements represented by some loops $\mu_i$ around $E_i$.

\begin{thm}\label{Mum} {\rm (\cite{Mu})}
The group $\pi_1^{loc}(B,o)$ of a curve germ $B$ is generated by $n+m$ elements $b_1,\dots, b_m$ and {\rm  (renumbered)} $e_{1+m},\dots, e_{n+m}$  being in one-to-one correspondence with the vertices of $\Gamma(B)$ and  being subject to the relations:
$$ \begin{array}{rll}
e_{i+m}^{w_{i+m}}\cdot b_1^{\delta_{1,i+m}}\cdot .\, .\, .
\cdot b_m^{\delta_{m,i+m}}\cdot
e_{1+m}^{\delta_{i+m,1+m}}\cdot . \, .\, .\cdot
e_{n+m}^{\delta_{i+m,n+m}}
&  =1,\quad & \text{for}\,\,i =1,\dots, n, \\
{[}b_j, e_{i+m}{]} & =1 \quad & \text{if}\,\,\,\, \delta_{j,i+m}=1,     \\
{[}e_{i_1+m}, e_{i_2+m}{]} & =1 \quad & \text{if}\,\,\,\, \delta_{i_1+m,i_2+m}=1.
\end{array}
$$
\end{thm}

\proof In \cite{Mu}, Mumford proved similar statement for presentations of the local fundamental groups of the complements to isolated two-dimensional singularities given in terms of resolution of singular points. In order to identify different fundamental groups, he defines a system of base paths lying in the curves $E_i$. We also  choose a system of paths on each $E_i$ as follows. Let
$$P_i=E_i\cap (\bigcup_{j\neq i}E_j)=\{ p_{i,j_1}=E_i\cap E_{j_1},\dots, p_{i,j_{k_i}}=E_i\cap E_{j_{k_i}}\},\,\, j_1<\dots <j_{k_i},$$
and
$$P'_i=E_i\cap (\bigcup B'_{j'})=\{ p_{i,j'_1}=E_i\cap B'_{j'_1},\dots, p_{i,j'_{k'_i}}=E_i\cap B'_{j'_{k'_i}}\},\,\, j'_1<\dots <j'_{k'_i}.$$
For each $E_i$, a point $Q_i\in E_i\setminus (P'_i\cup P_i)$ is connected by  paths $l'_{i,j'}$ with the points $p_{i,j'}\in P'_i$ and it is connected by  paths $l_{i,j}$ with the points $p_{i,j}\in P_i$. The paths $l_{i,j}$ and $l'_{i,j'}$ have the unique common point, namely $Q_i$, and if we move in a counter-clockwise direction along a circle $\delta\subset E_i$ of small radius with center at $Q_i$ 
then we consistently cross the paths $l'_{i,j'_1},\dots ,l'_{i,j'_{k'_i}}$ and then the paths $l_{i,j_1},\dots ,l_{i,j_{k_i}}$. After that everyone can easily check that the proof in \cite{Mu} can be transferred almost verbatim to the proof of Theorem \ref{Mum}.
\qed

\begin{prop} \label{A2n+1} If a curve germ $B$ has singularity type $A_{2n+1}$, $n\geq 0$, then 
$$ \pi_1^{loc}(B,o)=\langle b_1,b_2,e_3 \mid e_3=b_1b_2, [b_1,e_3^{n+1}]=[b_2,e_3^{n+1}]=1\rangle $$
and $e_{2+i}=e_3^i\,\, \text{for}\,\, i=2,\dots n+1.$ In particular, $e_{n+3}$ belongs to the center of $ \pi_1^{loc}(B)$.
\end{prop}
\proof  By Theorem \ref{Mum}, the group $\pi_1^{loc}(B,o)$ is generated by $b_1,b_2,e_3,\dots, e_{n+3}$ (see Fig. 1) and
\begin{equation} \label{A2n+1,1}
 e_3^2=e_4,\,\, e_4^2=e_3e_5,\,\, \dots, \,\, e_{n+2}^2=e_{n+1}e_{n+3},\end{equation}
\begin{equation} \label{A2n+1,2}  e_{n+3}=b_1b_2e_{n+2}, \end{equation}
\begin{equation}\label{A2n+1,3} {[e_{2+i}, e_{3+i}]}=1 \,\, \text{for} \,\, i= 1,\dots, n, \end{equation}
\begin{equation}\label{A2n+1,4} {[e_{n+3},b_{1}]}={[e_{n+3},b_2]}=1.
\end{equation}
It follows from (\ref{A2n+1,1}) that
$$e_4=e_3^2,\,\,\, \dots, \,\,\, e_{n+2}=e_3^n,\,\,\, e_{n+3}=e_{3}^{n+1}.$$
Therefore, by (\ref{A2n+1,2}) and (\ref{A2n+1,4}), $b_1b_2=e_3^{n+1}$ and ${[e_{n+3},b_{1}]}={[e_{n+3},b_2]}=1$. \qed

\begin{prop} \label{A2n} If a curve germ $B$ has singularity type $A_{2n}$, $n\geq 1$, then 
$$ \pi_1^{loc}(B,o)=\langle b_1,e_2 \mid e_2^{n+1}=b_1e_2^nb_1, [b_1,e_2^{2n+1}]=1\rangle $$
and $e_{1+i}=e_2^i\,\, \text{for}\,\, i=2,\dots n$, $e_{n+2}=e_2^{2n+1}$, and $e_{n+3}=b_1e_2^n$. In particular, $e_{n+2}$ belongs to the center of $ \pi_1^{loc}(B)$.
\end{prop}
\proof  By Theorem \ref{Mum}, the group $\pi_1^{loc}(B,o)$ is generated by $b_1,e_2,\dots, e_{n+3}$ (see Fig. 2) and
\begin{equation} \label{A2n,1}
 e_2^2=e_3,\,\, e_3^2=e_2e_4,\,\, \dots, \,\, e_{n}^2=e_{n-1}e_{n+1},\,\, e_{n+1}^3=e_ne_{n+2}\end{equation}
\begin{equation} \label{A2n,2}  e_{n+2}=b_1e_{n+1}e_{n+3}, \end{equation}
\begin{equation} \label{A2n,5}  e_{n+3}^2=e_{n+2}, \end{equation}
\begin{equation}\label{A2n,3} {[e_{1+i}, e_{2+i}]}=1 \,\, \text{for} \,\, i= 1,\dots, n, \end{equation}
\begin{equation}\label{A2n,4} {[e_{n+2},b_{1}]}={[e_{n+2},e_{n+1}]}={[e_{n+2},e_{n+3}]}=1.
\end{equation}
It follows from (\ref{A2n,1}) that
$$e_3=e_2^2,\,\,\, \dots, \,\,\, e_{n+1}=e_2^n,\,\,\, e_{n+2}=e_{2}^{2n+1}.$$
In particular, ${[e_{n+2},e_2]}=1$. Therefore, by (\ref{A2n,2}) and (\ref{A2n,5}), $e_{n+3}=b_1e_{n+1}^n$ and, by(\ref{A2n,4}), ${[e_{n+2},b_{1}]}={[e_{n+2},e_{n+3}]}=1$. \qed

\begin{prop} \label{D2n+2} If a curve germ $B$ has singularity type $D_{2n+2}$, $n\geq 1$, then 
$$ \pi_1^{loc}(B,o)=\langle b_1,b_2,b_3 \mid \,\,\, [b_1,b_2b_3]= [b_2,b_1(b_2b_3)^{n}]= [b_3,b_1(b_2b_3)^{n}]=1\rangle   $$
and $e_{3+i}=b_1(b_2b_3)^{i}$ for $i=1,\dots n$.
In particular, $e_{n+3}$ belongs to the center of $ \pi_1^{loc}(B)$.
\end{prop}
\proof  By Theorem \ref{Mum}, the group $\pi_1^{loc}(B,o)$ is generated by $b_1,b_2,b_3,e_4,\dots, e_{n+3}$ (see Fig. 3) and
\begin{equation} \label{D2n+2,1}
 e_4^2=b_1e_5,\,\, e_5^2=e_4e_6,\,\, \dots, \,\, e_{n+2}^2=e_{n+1}e_{n+3},\end{equation}
\begin{equation} \label{D2n+2,2}  e_{n+3}=b_2b_3e_{n+2}, \end{equation}
\begin{equation} \label{D2n+2,5}  {[b_1,e_{4}]}=1, \end{equation}
\begin{equation}\label{D2n+2,3} {[e_{2+i}, e_{3+i}]}=1 \,\, \text{for} \,\, i= 1,\dots, n, \end{equation}
\begin{equation}\label{D2n+2,4} {[e_{n+3},b_{2}]}={[e_{n+3},b_3]}={[e_{n+3},e_{n+2}]} =1.
\end{equation}
It follows from (\ref{D2n+2,5}) and (\ref{D2n+2,1}) that
\begin{equation}\label{D2n+2,6} e_5=b_1^{-1}e_4^2,\,\,\, \dots, \,\,\, e_{n+2}=b_1^{2-n}e_4^{n-1},\,\,\, e_{n+3}=b_1^{1-n}e_{4}^{n}.\end{equation}
Therefore, by (\ref{D2n+2,5}) and(\ref{D2n+2,6}), ${[e_{3+i},b_1]}=[e_{3+i},e_4]=1$ for $i=1,\dots,n$.
 It follows from (\ref{D2n+2,2}), (\ref{D2n+2,5}), and equalities $e_{n+2}=b_1^{2-n}e_4^{n-1},\,\,\, e_{n+3}=b_1^{1-n}e_{4}^{n}$ that $e_4=b_1b_2b_3$ and $[b_1,b_2b_3]=1$. By (\ref{D2n+2,6}), we have $e_{3+i}=b_1(b_2b_3)^i$ for $i=1,\dots, n$. In particular, $e_{n+3}=b_1(b_2b_3)^{n}$ and, by (\ref{D2n+2,4}), $e_{n+3}$ belongs to the center of $ \pi_1^{loc}(B)$. Now, it follows from (\ref{D2n+2,5}) that $[b_1,b_2b_3]=1$ and it is easy to see that relations (\ref{D2n+2,3}) do not give additional relations. \qed \\

The proofs of the following four Lemmas are similar to the proofs of Lemmas \ref{A2n+1} -- \ref{D2n+2} and therefore they they will be omitted.

\begin{prop} \label{D2n+3} If a curve germ $B$ has singularity type $D_{2n+3}$, $n\geq 1$, then {\rm (see Fig. 4)}
$$\pi_1^{loc}(B,o)=\langle b_1,b_2,e_3 \mid
e_3^{2n}b_1^{1-2n}=(e_3^nb_1^{-n}b_2^{-1})^2,\,\, [e_3,b_1]= [b_2,e_3^{2n}b_1^{1-2n}]=1   \rangle
$$
and $e_{i+2}=e_2^{i}b_1^{1-i}$ for $i=1,\dots n$, $e_{n+3}=e_3^{2n}b_1^{1-2n}$, and $e_{n+4}=e_3^nb_1^{-n}b_2^{-1}$.
In particular, $e_{n+3}$ belongs to the center of $ \pi_1^{loc}(B)$.
\end{prop}

\begin{prop} \label{E6} If a curve germ $B$ has singularity type $E_{6}$,  then {\rm (see Fig. 5)}
$$
\pi_1^{loc}(B,o)=\langle b_1,e_2 \mid 
e_2^{3}=(b_1e_2)^2b_1,\,\, [e_2^4,b_1]=1
\rangle
$$
and $e_{3}=e_2^{4}$, $e_{4}=(b_1e_2)^2$, and $e_{5}=b_1e_2$. In particular, $e_{3}$ belongs to the center of $ \pi_1^{loc}(B)$.
\end{prop}

\begin{prop} \label{E7} If a curve germ $B$ has singularity type $E_{7}$,  then {\rm (see Fig. 6)}
$$
\pi_1^{loc}(B,o)=\langle b_1,b_2,e_3 \mid 
e_3^{2}=b_1b_2e_3b_2,\,\, [b_1,b_2e_3]=[e_3^3,b_1]=[e_3^3,b_2]
\rangle
$$
and $e_{4}=e_3^{3}$, $e_{5}=b_1b_2e_3$. In particular, $e_{4}$ belongs to the center of $ \pi_1^{loc}(B)$.
\end{prop}

\begin{prop} \label{E8} If a curve germ $B$ has singularity type $E_{8}$,  then {\rm (see Fig. 7)}
$$
\pi_1^{loc}(B,o)=\langle b_1,e_2 \mid (e_2^{2}b_1^{-1})^2=b_1e_2^3,\,\, [e_2^5,b_1]=1 \rangle
$$
and $e_{3}=e_2^{3}$, $e_{4}=e_2^{5}$, $e_5=e_2^2b_1^{-1}$.
In particular, $e_{4}$ belongs to the center of $ \pi_1^{loc}(B)$.
\end{prop}

\begin{cor} \label{corol} Let $(B,o)$ be a curve germ having one of $ADE$ singularity types, $E\subset \sigma^{-1}(o)\subset V_n$ the exceptional curve of the last blowup $\sigma_n$ 
in the sequence of blowups resolving the singular point of $(B,o)$, and $e$ an  element in $\pi_1^{loc}(B,o)$ represented by a simple loop around $E$. Then $e$ belongs to the center of $\pi_1^{loc}(B,o)$.
\end{cor}
\begin{prop} \label{gener} Let $(B,o)$ be a curve germ having one of $ADE$ singularity types.  If the singularity type of $(B,o)$ is not $A_0$ or $A_1$, then $\pi_1^{loc}(B,o)$ is generated by  $e$ and the elements $\gamma_1,\gamma_2,\gamma_3$ corresponding to the vertices of $\Gamma(B)$ connected by an edge with the vertex $e$ {\rm (if the singularity type of $(B,o)$ is  $A_1$, then $\pi_1^{loc}(B,o)$ is generated by  $b_1,b_2$ and $e$)}.
\end{prop}
\proof We prove Proposition \ref{gener} only in the case when the singularity type of $(B,o)$ is $D_{2n+3}$.  Proof in all other cases is similar and will be omitted.

In the case of singularity type $D_{2n+3}$, the elements $\gamma_1,\gamma_2,\gamma_3$ are $b_2,e_{n+2}, e_{n+4}$ and $e=e_{n+3}$ (see Fig. 4). By Theorem \ref{Mum}, the group $\pi_1^{loc}(B,o)$ is generated by the elements $b_1,b_2, e_3,\dots, e_{n+4}$ and among the relations connecting these elements, we have the following relations:
$$ e_3^2=b_1e_4,\,\, e_4^2=e_3e_5,\,\,\dots,\,\, e_{n+1}^2=e_ne_{n+2},\,\, e_{n+2}^3=e_{n+1}e_{n+3}.$$ Therefore
$e_{n+2-i}=e_{n+2}^{2i+1}e_{n+3}^{-1}$ for $i=1,\dots, n-1$ and $b_1=e_{n+2}^{2n+1}b_{n+3}^{-1}$. \qed \\

Denote by $Z_e$ the subgroup of $\pi_1^{loc}(B,o)$ generated by $e$ and consider the group
$\pi_1(V_n\setminus \overline{\sigma^{-1}(B)\setminus E}),p)$, where $\overline{\sigma^{-1}(o)\setminus E}$ is the closure of $\sigma^{-1}(o)\setminus E$ in $V_n$. Without loss of generality, we can assume that $p\in E$. Corollary \ref{corol} and Proposition \ref{gener} imply
\begin{claim} \label{epi} A homomorphism
$$i_*: \pi_1(E\setminus \overline{\sigma^{-1}(B)\setminus E},p)\to \pi_1(V_n\setminus \overline{\sigma^{-1}(B)\setminus E},p)\simeq \pi_1^{loc}(B,o)/Z_e,$$
induced by imbedding
$i: E\setminus \overline{\sigma^{-1}(B)\setminus E}\hookrightarrow V_n\setminus \overline{\sigma^{-1}(B)\setminus E}$, is an epimorphism.
\end{claim}

\section{Proof of Theorem 2}
Consider a finite cover $F:(U,o')\to (V,o)$, $F\in \mathcal R$, and let $(B,o)$ be its branch curve, $G_F\subset \mathbb S_d$ its monodromy group, where $d=\deg_{o'} F$. Let $\sigma: V_n\to V$ be the minimal resolution of the singular point $o$ of $B$. 
The  cover $F:U\setminus F^{-1}(B)\to V\setminus B$ is unramified and  by Grauert - Remmert - Riemann - Stein Theorem, the monodromy $$F_*:\pi_1^{loc}(B,o)\simeq \pi_1(V\setminus B,p)\simeq \pi_1(V_n\setminus \sigma^{-1}(o),\sigma^{-1}(p))\to \mathbb S_d$$  defines a finite holomorphic map $F_n:U_n\to V_n$ branched in $\sigma^{-1}(B)$, where $U_n$ is a normal complex-analytic variety such that if $U_n$ is singular, then its singular points lie over the singular points of the divisor $\sigma^{-1}(B)$ (and they are singularities of Hirzebruch-Jung type (see \cite{B})). In addition, the cover
$F:U\setminus F^{-1}(B)\to V\setminus B$ and the cover $F_n:U_n\setminus F_n^{-1}(\sigma^{-1}(B))\to V_n\setminus \sigma^{-1}(B)$  are the same cover. Therefore, by Stein factorization theorem applied to the map $\sigma\circ F_n$, there is
a commutative diagram
\begin{equation} \label{diagram} \begin{CD} U_n @>{F_n}>> V_n \\ @V{\Sigma}VV @VV{\sigma}V \\
U @>{F}>> V
\end{CD}\end{equation}
in which
$\Sigma: U_n\to U$ is  a holomorphic bimeromorphic map contracting the curves lying in $F_n^{-1}(\sigma^{-1}(o))$ to the point $o'$.
By Zariski Theorem applied to the composition of the minimal resolution of singular points of $U_n$ and $\Sigma$, all curves from $F_n^{-1}(\sigma^{-1}(o))$ are rational. Note also that $F_n^{-1}(\sigma^{-1}(o))=\Sigma^{-1}(o')$ is connected.

Let $E\subset \sigma^{-1}(o)\subset V_n$ be the exceptional curve of the last blowup $\sigma_n$  and $e$ an  element in $\pi_1^{loc}(B,o)$ represented by a simple loop around $E$. Then, by Corollary \ref{corol}, $e$ belongs to the center of $\pi_1^{loc}(B,o)$. Denote by
\begin{equation} \label{Z} Z:=F_*(Z_e)\subset G_F.\end{equation}
where $Z_e$ is the subgroup of $\pi_1^{loc}(B,o)$ generated by $e$. Note that $Z\neq G_F$ for $F\in\mathcal R\setminus (\mathcal R_{A_0}\cup\mathcal R_{A_1})$.

\begin{prop} \label{factor} Let  $F:X\to Y$, $\deg F=d$, be a  finite holomorphic map from a connected normal complex-analytic variety $X$ to a smooth complex surface $Y$ branched in a curve $B\subset  Y$ and $Z$ a subgroup of the center of the monodromy group $G_F\subset \mathbb S_d$ of $F$. Then
\begin{itemize}
\item[($i$)] the order $|Z|$ of $Z$ is a divisor of $\deg F$, $d=d_1\cdot |Z|$,
\item[($ii$)] $Z$ acts on $X$ and the quotient variety $W=X/Z$ is  a normal 
variety,
\item[($iii$)] $F=H\circ F_Z$, where $F_Z:X\to W$ is the quotient map, $\deg F_Z=|Z|$, and $H:W\to Y$ is a holomorphic finite map,
$\deg H=d_1$.
\end{itemize}
\end{prop}
\proof
Consider the symmetric group $\mathbb S_d$ as a group acting on the interval of natural numbers $\mathbb N_d=\{ 1,\dots, d\}$. Denote by $\mathbb S_{d-1}$ the subgroup of $\mathbb S_d$ consisting of the permutations $\tau\in\mathbb S_d$ leaving fixed $1$. Since $G_F$ is a transitive subgroup of $\mathbb S_d$, then $G_1=G_F\cap \mathbb S_{d-1}$ is a subgroup of $G_F$ of index $(G_F:G_1)=d$. Let us show that $G_1$ is a relatively simple subgroup of $G_F$, i.e. $G_1$ does not contain a proper non-trivial normal subgroup of $G_F$. Indeed, assume that a normal subgroup $N$ of $G_F$ is contained in $G_1$ and $h\in N$ is a non-trivial element. But,  for any $i\in \mathbb N_d$ there is an element $g_i\in G_F$ such that $g_i(1)=i$ and therefore $h(i)=i$ for each $i\in\mathbb N_d$, since $g_i^{-1}hg_i\in N\subset G_1$. As a result, we get a contradiction with the assumption that $G_F\subset \mathbb S_d$.

Let $c:G=G_F\hookrightarrow \mathcal S_{|G|}$ be  Cayley's imbedding. By  Grauert - Remmert - Riemann - Stein Theorem, the homomorphism $c\circ F_*: \pi_1(Y\setminus B,p)\to \mathbb S_{|G|}$ defines a {\it  Galois cover} $\widetilde F: \widetilde X\to Y$ of degree
$\deg \widetilde F=|G|$, where $\widetilde F$ is a holomorphic finite map and  $\widetilde X$ is  a connected normal complex-analytic variety. The group $G$ acts on $\widetilde X$ such that the quotient variety
$\widetilde X/G$ is $Y$ and $\widetilde F$ is the quotient map. It is well known that  the quotient variety
$\widetilde X/G_1$ is biholomorphic to $X$ and the map $\widetilde F$ is the composition of two maps, $\widetilde F=F\circ  F_{G_1}$, where $F_{G_1}: \widetilde X\to X$ is the quotient map defined by the action of $G_1$ on $\widetilde X$.

Denote by $\widetilde G_1=G_1Z$ the subgroup of $G$ generated by the elements of $G_1$ and $Z$. Then $\widetilde F=H\circ F_{\widetilde G_1}$, where $F_{\widetilde G_1}:\widetilde X\to W$ is the quotient map, $\deg F_{\widetilde G_1}=|\widetilde G_1|$, and $H:W\to Y$ is a holomorphic finite map.

Since $Z$ is a normal central subgroup and $G_1$ is a relatively simple subgroup of $G$, then $G_1\cap Z=\{ {\bf 1}\}$. Therefore $\widetilde G_1$ is isomorphic to $G_1\times Z$, the group $G_1$ is a normal subgroup of $\widetilde G_1$, and hence $F_{\widetilde G_1}=F_Z\circ F_{G_1}$, where
$F_Z:X\to W$ is the quotient map defined by the action of the group $Z=\widetilde G_1/G_1$ on $X$, $\deg F_{Z}=|Z|$. Now, Proposition \ref{factor} follows from the equalities $\widetilde F=H\circ F_{\widetilde G_1}=H\circ F_Z\circ F_{G_1}$ and $\widetilde F=F\circ F_{G_1}$. \qed

\begin{rem} \label{rem2} The monodromy group of the finite cover $H$ in Proposition {\rm \ref{factor}} is $G_{H}=G_F/N\subset \mathbb S_{d_1}$, where $N$ is the maximal normal subgroup of $G_F$ contained in $G_1Z$ {\rm (the group $G_1$ is defined in the proof of Proposition \ref{factor})}.
\end{rem}

Let us return to the case when $F\in \mathcal R$ and apply Proposition \ref{factor} to diagram (\ref{diagram}) (the cyclic group $Z\subset G_F$ 
is defined in (\ref{Z})). As a result, we obtain the following commutative diagram:
$$\begin{CD} F_n:\phantom{a}
@. U_n @>{F_{n,Z}}>>W_n @>H_n>>  V_n @. \phantom{a}\supset E \\ @.  @V{\Sigma}VV @V{\Theta}VV @VV{\sigma}V \\
F:
@. U @>{F_Z}>> W @>{H}>> V @. \ni o
\end{CD}
$$
in which $H$ and $H_n$ are finite holomorphic maps, $\deg H_n=\deg H =d_1=\frac{d}{|Z|}$, and $\Theta$ contacts $H_n^{-1}(\sigma^{-1}(o))$ to the point $o_1=H^{-1}(o)=F_Z(o')$. The proper inverse image $H_n^{-1}(\sigma^{-1}(o))$ is a union of rational curves, since $F_n^{-1}(\sigma^{-1}(o))$ is a union of rational curves. The monodromy group $G_{H_n}$ of $H_n$ is $G_F/N\subset \mathbb S_{d_1}$, where $N$ is a normal subgroup of $G_F$ defined in Remark \ref{rem2}, and the monodromy homomorphism $H_{n*}:\pi_1(V_n\setminus \sigma^{-1}(B),p)\to G_F/N$ is a composition of the following homomorphisms: $F_*:\pi_1(V_n\setminus \sigma^{-1}(B),p)\to G_F$, the quotient homomorphism  $G_F\to G_F/N$, and an embedding $G_F/N\hookrightarrow \mathbb S_{d_1}$.  The map $H_n$ is not branched in $E$, since $F_*(e)\in Z\subset N$. Therefore $H_{n*}$ can be considered as a homomorphism
$$H_{n*}:\pi_1((V_n\setminus \sigma^{-1}(B))\cup E,p)\to G_F/N.$$

The intersection matrix of the irreducible components of the closure $\overline{\sigma^{-1}(o)\setminus E}$ of $\sigma^{-1}(o)\setminus E$ in $V_n$  is negatively defined. Therefore  $\sigma=\varphi\circ\psi$, where  $\psi: V_n\to S$ is the contraction contracting the divisor $\overline{\sigma^{-1}(o)\setminus E}$ to points and  $\varphi:S\to V$ is the holomorphic map  contracting $\psi(E)$ to the point $o$. Note that $\psi_{\mid E}:E\to \psi(E)$ is an isomorphism.

By Stein factorization theorem, $\psi\circ H_n=\beta\circ \xi$, where  $\xi: W_n\to T$ is the contraction contracting the divisor $H_n^{-1}(\overline{\sigma^{-1}(o)\setminus E})$ to points and   $\beta:T\to S$ is a finite holomorphic map, $\deg \beta=d_1$ and the monodromy group $G_{\beta}=G_F/N$.

It is easy to see that $\xi_{\mid H_n^{-1}(E)}:H_n^{-1}(E)\to \xi(H_n^{-1}(E))$ is an isomorphism. Therefore $\xi(H_n^{-1}(E))=\mathbb P^1$, since $H_n^{-1}(\sigma^{-1}(o))$ is a connected union of rational curves, We obtain a finite holomorphic map
$f=\beta_{\mid \xi(H_n^{-1}(E))}:\xi(H_n^{-1}(E))\simeq \mathbb P^1\to \psi(E)\simeq \mathbb P^1$
branched in no more than  three points $\psi(\overline{\sigma^{-1}(o)\setminus E})\subset \psi(E)$, $\deg f =d_1$.
\begin{df} \label{Be}  The map $\beta:\mathcal R \to \mathcal Bel$ sends $F\in\mathcal R$ to $\beta(F)\in \mathcal Bel$ by the following rule:
\newline $\bullet$ if $F\in \mathcal R_{A_0}$, then $\beta(F)=\text{id}:\mathbb P^1\to\mathbb P^1\in\mathcal Bel_{2}$;
\newline $\bullet$ if $F\in \mathcal R_{A_1}$ with $G_F=\mathbb Z_{n_1}\times\mathbb Z_{n_2}$, $\text{GCD}(n_1,n_2)=k$, then $\beta(F)\in \mathcal Bel_{2}$ with $\phantom{aa}G_{\beta(F)}=\mathbb Z_k$;
\newline $\bullet$ if $F\in\mathcal R\setminus (\mathcal R_{A_0}\cup\mathcal R_{A_1})$, then 
$$\beta(F):=\beta_{\mid \xi(H_n^{-1}(E))}:\xi(H_n^{-1}(E))\simeq \mathbb P^1\to \psi(E)\simeq \mathbb P^1.$$
\end{df}
It easily follows from
Claim \ref{epi} that the monodromy group $G_{\beta(F)}\simeq G_{H_n}\simeq G_F/N$ if $F\in\mathcal R\setminus (\mathcal R_{A_0}\cup\mathcal R_{A_1})$.

To complete proof of Theorem 2, let us show that for each $f\in\
mathcal Bel_3$ of degree $\deg f=n$, there is a finite cover $F\in \beta^{-1}(f)\cap\mathcal R_{D_4}$ of degree $\deg_o F=n^2$.

Let $$T_c(f)=\{ c_i=(m_{1,i},\dots ,m_{k_i,i})\}_{m_{1,i}+\dots + m_{k_i,i}=\deg f,\,\, i\in \{ 0,1,\infty\}}$$ be the the set of cycle types of permutations $f_*(\gamma_i)\in G_f\subset \mathbb S_n$,  $i\in \{ 0,1,\infty\}$ and let $(B,o)$ has the singularity type $D_4$. Its local fundamental group is described in Proposition \ref{D2n+2}. Consider a homomorphism $H_{1*}:\pi_1^{loc}(B,o)\simeq \pi_1(V_1\setminus \sigma_1^{-1}(B))\to \mathbb S_n$ sending $b_1$ to $f_*(\gamma_0)$, $b_2$ to $f_*(\gamma_1)$, and $b_3$ to $f_*(\gamma_{\infty}^{-1})$.
We have $H_{1*}(e)=H_{1*}(b_1b_2b_3)=\text{id}$, where $\text{id}\in \mathbb S_n$ is the identical permutation. The homomorphism $H_{1*}$ defines a finite covering $H_1:W_1\to V_1$ branched in $B'_1\cup B'_2\cup B'_3$ and it does not ramified over $E$, since $H_{1*}(e)=\text{id}$. Therefore $W_1$ is a smooth surface, since $H_1$ is branched in the disjoint union of smooth curve germs.

By  Claim \ref{epi},  $\widetilde E=H_1^{-1}(E)\simeq \mathbb P^1$ and  $H_{1\mid \widetilde E}=f$. The intersection number $(\widetilde E^2)_{W_1}=\deg H_1\cdot (E^2)_{V_1}=-n$. Therefore $\pi_1(W_1\setminus \widetilde E)=Z\simeq\mathbb Z_n$ and hence there is a cyclic cover
$F_{1,Z}:U_1\to W_1$ branched in $\widetilde E$ with multiplicity $n$. Let $\overline E$ be the proper inverse image $F_{1,Z}^{-1}(\widetilde E)\simeq \widetilde E\simeq \mathbb P^1$. We have
$$(F_{1,Z}^*(\widetilde E),F_{1,Z}^*(\widetilde E))_{U_1}=(n\overline E,n\overline E)_{U_1}=\deg F_{1,Z}\cdot (\widetilde E^2)_{W_1}=-n^2.$$
Therefore $(\overline E^2)_{U_1}=-1$ and hence there is the contraction ($\sigma$-process) $\Sigma_1: U_1\to U$ contracting $\overline E\simeq \mathbb P^1$ to a smooth point. It is easy to see that $$F:=\sigma_1\circ H_1\circ F_{1,Z}\circ \Sigma_1^{-1}:U\to V$$ is a finite cover, $\deg F=n^2$. The cover $F$  is branched in $(B,o)$ and $\beta(F)=f$.

\section{Proof of Theorem \ref{main3}}

Without loss of generality, we can assume that the branch curve germ  $(B,o)\subset (V,o)\subset (\mathbb C^2,o)$  of the cover
$F:(U,o')\to (V,o)$, $\deg_{o'} F=d$, is  given in $\mathbb C^2$ by one of the following equations:
\begin{itemize}
\item[($A_n$)] \phantom{aaa} $u^2-v^{n+1}=0$, $n\geq 0$;
\item[($D_n$)] \phantom{aaa} $v(u^2-v^{n-2})=0$, $n\geq 4$;
\item[($E_6$)] \phantom{aaa} $u^3-v^4=0$;
\item[($E_7$)] \phantom{aaa} $u(u^2-v^3)=0$;
\item[($E_8$)] \phantom{aaa} $u^3-v^5=0$
\end{itemize}
and $(B,o)$ is the germ at $o=(0,0)\in \mathbb C^2$ of an affine curve $B\subset \mathbb C^2$ given in coordinates $(u,v)$ by the same equation. Let $(z_0:z_1:z_2)$ be homogeneous coordinates in $\mathbb P^2$ and $\mathbb C^2\hookrightarrow\mathbb P^2$ an imbedding given by $u=\frac{z_0}{z_2}$, $v\frac{z_1}{z_2}$. Denote by  $\overline B$ the closure of $B$ in $\mathbb P^2$ and $L_i\subset\mathbb P^2$, $i=0,1,2$, a line given by  equation $z_i=0$.

Note that the equations ($A_n$) -- ($E_8$) are quasi-homogeneous. Therefore
$$\pi_1^{lpc}(B,o)\simeq \pi_1(\mathbb C^2\setminus B)=\pi_1(\mathbb P^2\setminus (\overline B\cup L_{2})).$$
and the monodromy homomorphism $F_*:\pi_1^{loc}(B,o)\simeq \pi_1(\mathbb P^2\setminus (\overline B\cup L_{2}))\to \mathbb S_d$ defines a finite
cover $\overline F:X\to \mathbb P^2$ branched in $\overline B\cup L_{2}$, $\deg \overline F=d$, where $X$ is  a normal irreducible complex-analytic surface.  Obviously, $F:(U,o')\to (V,o)$ is a germ of the cover $\overline F$, $(U,o')=(\overline F^{-1}(V),\overline F^{-1}(o))\subset (X,o')$. Therefore $X$ is smooth at $o'$.

Obviously, Theorem \ref{main3} holds if $F\in \mathcal R_{A_0}\cup\mathcal R_{A_1}$. Therefore we will assume that $F\in \mathcal R\setminus (\mathcal R_{A_0}\cup\mathcal R_{A_1})$.

Consider, in the beginning, the case when $F\in \mathcal R_{D_4}$. Denote by $L_{0,1}\subset \mathbb P^2$ the line given  by equation $z_0-z_1=0$, then $(B,o)$ is the germ of curve $\overline B=L_{0}\cup L_{1}\cup L_{0,1}$ at $o=(0,0,1)$. Let $\sigma=\sigma_0\circ\sigma_{1}:Y_1\to \mathbb P^2$ be two $\sigma$-processes with centers at $o=(0,0,1)$ and $o_1=(1,1,0)$. Denote by the same letters the proper inverse images of the lines $L_{0}$,
$L_{1}$,   $L_{2}$, $L_{0,1}$ and let $E=\sigma^{-1}(o)$, $E_1=\sigma^{-1}(o_1)$. After that we blowdown the curve $L_{0,1}$ to a point by $\sigma$-process $\tau:Y_1\to Y_2$. It is easy to see that $Y_2$ is isomorphic to $\mathbb P^1\times\mathbb P^1$ in which $L_{0}$, $L_{1}$, $E_1$ are fibres of the projection to the first factor and  $L_{2}$, $E$ are sections. Note that $Y_1$ and $Y_2$ are defined over
$\mathbb Q$.

We have
\begin{equation}\label{iso}\mathbb C^2\setminus B=Y_1\setminus (L_{0}\cup L_{1}\cup L_{2}\cup L_{0,1}\cup E\cup E_1)=Y_2\setminus (L_{0}\cup L_{1}\cup L_{2}\cup E\cup E_1)\end{equation}
(below, we identify these surfaces). Therefore
$$\pi_1(\mathbb C^2\setminus B)=\pi_1(Y_1\setminus (L_{0}\cup L_{1}\cup L_{2}\cup L_{0,1}\cup E\cup E_1))=\pi_1(Y_2\setminus (L_{0}\cup L_{1}\cup L_{2}\cup E\cup E_1)).$$
The monodromy homomorphism $F_*:\pi_1^{loc}(B,o)\to \mathbb S_d$ defines finite covers $F:X\to \mathbb P^2$ unramified over $\mathbb C^2\setminus B\subset \mathbb P^2$ and $F_i:X_i\to Y_i$, $i=1,2$, unramified over $\mathbb C^2\setminus B\subset Y_i$. Denote by $\nu:\widetilde X\to X$ and $\nu_i:\widetilde X_i\to X_i$, $i=1,2$, resolutions of singular points of $X$ and $X_i$.
We have the following commutative diagram
\begin{equation} \label{diag2} \begin{CD}
\widetilde X @<\widetilde{\sigma}<< \widetilde X_1 @>\widetilde{\tau}>> \widetilde X_2  \\   @V{\nu}VV @V{\nu_1}VV @VV{\nu_2}V \\
 X @<\overline{\sigma}<<X_1 @>\overline{\tau}>>  X_2  \\   @V{F}VV @V{F_1}VV @VV{F_2}V \\
\mathbb P^2 @<{\sigma}<< Y_1 @>{\tau}>> Y_2 @. =\mathbb P^1\times\mathbb P^1
\end{CD}
\end{equation}
in which all horizontal arrows are bimeromorphic maps and the finite covers $F$, $F_1$, $F_2$ are
the same unramified cover over $\mathbb C^2\setminus B$. Since $o'$ is a smooth point of $X$, then $\nu :\nu^{-1}(U)\to U$ is a biholomorphic map and therefore we will identify  $(U,o')$ with $(\nu^{-1}(U),\nu^{-1}(o'))$  and $F:(U,o')\to (V,o)$ with the restriction to $(\nu^{-1}(U),\nu^{-1}(o'))$ of the map $F\circ\nu: \widetilde X\to \mathbb P^2$.

It is easy to see that all surfaces  included in diagram (\ref{diag2}) are algebraic surfaces and all maps between them are regular morphisms. Indeed,  $\widetilde X$ and $\widetilde X_i$, $i=1,2$, are projective surfaces, since the transcendence degrees of the fields of meromorphic functions $\mathbb C(\widetilde X)$ and $\mathbb C(\widetilde X_i)$ equal two (these fields contain the field $\mathbb C(\mathbb P^2)$), and the varieties $X$ and $X_i$ coincide, resp., with the normalizations of $\mathbb P^2$ and $Y_i$, $i=1,2$, in the fields  $\mathbb C(\widetilde X)\simeq\mathbb C(\widetilde X_i)$ (see Chapter II, Section 5, Subsection 2 in \cite{Sh}).

Let us show that $\widetilde X_2$ is a rational surface defined over $\overline{\mathbb Q}$. For this, consider again  the subgroup $Z_e$ of $\pi_1^{loc}(B,o)$ and its image $Z=F_*(Z_e)$ (see Corollary \ref{corol}). By Proposition \ref{factor}, we have $F_2=H\circ F_{2,Z}$, where $F_{2,Z}:X_2\to W$ is the quotient map under the action of cyclic group $Z$ on $X_2$, $\deg F_{2,Z}=|Z|$, and $H:W\to Y_2=\mathbb P^1\times\mathbb P^1$ is a finite morphism branched only in three fibres $L_0$, $L_{1}$, and $E_1$ (over the points $\{ 0,1,\infty\}\in L_2$) of the first projection. According to Definition \ref{Be}, by Theorem \ref{main2}, the restriction of $H$ to $H^{-1}(E)$ is a Belyi function $\beta(F):H^{-1}(E)\simeq\mathbb P^1\to\mathbb P^1$. 
The inverse image $H^{-1}(C)$ of a generic fiber $C$ of the first projection is the disjoint union $\bigsqcup_{i=1}^{\deg H}C_i$ of curves isomorphic to $C$ and $(H^{-1}(E),H^{-1}(E))_W=0$. Therefore $W\simeq H^{-1}(E)\times C_1\simeq \mathbb P^1\times \mathbb P^1$ and
$$H=\beta(F)\times\text{id}:W\simeq H^{-1}(E)\times C_1\to \mathbb P^1\times C
.$$
Hence, $H:W\to Y_2$ and $R=H^{-1}(L_0\cup L_1\cup E_1)$, $H^{-1}(E)$, $H^{-1}(L_2)$ are defined over $\overline{\mathbb Q}$.

The cyclic cover $\overline F_{2,Z}$, $\deg F_{2,Z}=|Z|=n$,  is branched with multiplicity $n$ in sections $H^{-1}(E)$, $H^{-1}(L_2)$ of the ruled structure on $W$ defined by the projection to the first factor, and, possibly, in several fibres belonging to $R$.  Therefore $\widetilde X_2$ can have only cyclic quotient singularities (locally the normalizations   of singularities given by $z_i^{k_i}=x_iy_i$, where $k_i$ are divisors of $n$, see \cite{B}) over points in $R\cap H^{-1}((E\cup L_{2})$, and the curves $F_{2,Z}^{-1}(H^{-1}(E))\simeq H^{-1}(E)$ and $F_{2,Z}^{-1}(L_1)$ are rational. As a result, we see that $\widetilde X_2$ has a ruled structure with a rational section. Therefore $\widetilde X_2$ is a rational surface defined over $\overline{\mathbb Q}$.

The surface 
$X_1$ is the normalization of the fibre product $Y_1\times_{Y_2} X_2$ of morphisms $\tau_1:Y_1\to Y_2$ and $F_2:\overline X_2\to Y_2$ defined over the field $\overline{\mathbb Q}$, and $X$ is the normalization of $\mathbb P^2$ in the field $\overline{\mathbb Q}(X_1)$ containing the field $(\tau\circ F_1)^*(\overline{\mathbb Q}(\mathbb P^2))$. Therefore the surfaces $X_1$, $X$,  and  their resolutions of singularities of these surfaces are rational and defined over $\overline{\mathbb Q}$.

\begin{claim}\label{rat} For each point $p$ of  a rational smooth projective surface $S$ defined over algebraicly closed field $\overline{\Bbbk}$, $\text{char}\, \overline{\Bbbk}=0$, 
there is a Zariski open neighbourhood $U\subset S$ of the point $p$ such that $U$ is isomorphic to the affine surface $\overline{\Bbbk}^2$.
\end{claim}
\proof Remind that if $x,y$ are coordinates in $\Bbbk^2$ and $\sigma: X\to \Bbbk^2$ is the $\sigma$-process with center at $o=(0,0)$, then
$X$ is covered by two Zariski open neighbourhoods $U_1$ and $U_2$ isomorphic to $\Bbbk^2$ and such that $\sigma:U_1\to\Bbbk^2$ in coordinates $x_1,y_1$ in $U_1$ is given by $x=x_1$ and $y=x_1y_1$, resp., $\sigma:U_2\to\Bbbk^2$ in coordinates $x_2,y_2$ in $U_2$ is given by $x=x_2y_2$ and $y=y_2$. Note also that for each point $p\in \Bbbk^2$ we can choose  coordinates $x,y$ in $\Bbbk^2$ such that $p$ is the origin of this coordinate system. To complete the proof of Claim \ref{rat} it suffices to remind that for each rational smooth projective surface $S$ there is a  birational morphism $f:S\to M$ to a relatively minimal model $M$ which is a composition of $\sigma$-processes by Zariski Theorem, where $M$ is isomorphic either to a Hirzebruch surface $\mathbb F_n$, or $\mathbb P^2$ (see \cite{Sh1}) and for each point $p\in M$ there is a Zariski open neighbourhood $U\subset M$ of $p$ such that $U$ is isomorphic to the affine surface $\overline{\Bbbk}^2$. \qed \\


By Claim \ref{rat}, we can choose a Zariski open neighbourhood $\widetilde U\subset \widetilde X$ defined over $\overline{\mathbb Q}$ isomorphic to $\mathbb C^2$ and containing the point $o'$. Choose coordinates $z,w\in \overline{\mathbb Q}[z,w]$ in $\widetilde U$ such that $o'=(0,0)$ is the origin of this coordinate system. Then the restriction of $F\circ\nu$ to $\widetilde U$ defines a rational map $\widetilde U\dashrightarrow \mathbb C^2\subset \mathbb P^2$ regular at $o'$. Therefore this rational map is given by functions
$$u=\frac{f_1(z,w)}{g_1(z,w)}, \quad v=\frac{f_2(z,w)}{g_2(z,w)},$$
where $f_i(z,w)$ and $g_i(z,w)\in \overline{\mathbb Q}[z,w]$ for $i=1,2$ and $g_1(0,0)g_2(0,0)\neq 0$. The restriction of this map to $(U,o')$ is $F:(U,o')\to (V,o)$.

To prove Theorem \ref{main3} in the case $F\in\mathcal R\setminus (\mathcal R_{A_0}\cup \mathcal R_{A_1}\cup\mathcal R_{D_4})$, consider in $\mathbb P^2$ two pencils of curves given in homogeneous coordinates $(z_0:z_1:z_2)$ in $\mathbb P^2$ by equations

\begin{equation}\label{P1} \lambda z_0z_2^{n-1}\,\,\,\,\,  +\,\,  \mu z_1^{n}\,\,  =  0,\qquad n\geq 2;\end{equation}
\begin{equation}\label{P2} \lambda z_0^2z_2^{2n-1}  +\,\,  \mu z_1^{2n+1}  =  0,\qquad n\geq 2.
\end{equation}
These pencils define two rational maps $\varphi_i:\mathbb P^2\to \mathbb P^1$, $i=1,2$. The maps $\varphi_i$ have two indeterminacy points $A_1=(0:0:1)$ and $A_2=(1:0:0)$. To resolve the indeterminacy points of pencil (\ref{P1}), we need to blow up $k=n$ times each point $A_1$ and  $A_2$ and to resolve the indeterminacy points of pencil (\ref{P2}), we need to blow up $k=n+2$ times each point $A_1$ and  $A_2$. Let $\sigma=\sigma_1\circ\dots\circ\sigma_k\circ\sigma_{k+1}\circ\dots\circ\sigma_{2k}:Y_1\to \mathbb P^2$ be a sequence of $\sigma$-processes resolving the indeterminacy points, where the first $k$ $\sigma$-processes blow up the point $A_1$ and the points lying over $A_1$. Denote by $E_j\subset Y$ the proper inverse image of the exceptional curve of $\sigma_j$, $j=1,\dots,2k$, and by $L_i\subset Y$, $i=0,1,2$, the proper inverse image of the line given in $\mathbb P^2$ by equation $z_i=0$.

It is easy to check that in  the case when the pencil is given by equation (\ref{P1}) and $n\geq 3$, the weighted dual graph of the irreducible components of the curve $C=(\bigcup E_j)\cup(\bigcup L_i)\subset Y$ is the graph depicted in Fig. 8, where the weights are the self-intersection numbers in $Y_1$ of the irreducible components of $C$, and in  the case when the pencil is given by equation (\ref{P1}) and $n=2$, the weighted dual graph of the curve $C=(\bigcup E_j)\cup(\bigcup L_i)\subset Y$  is the graph depicted in Fig. 9.

\begin{picture}(300,85)
\put(85,30){\circle*{3}}\put(70,30){$\mbox{}_{-1}$}
\put(88,25){$E_{2k}$}
\put(325,30){\circle*{3}}\put(330,30){$\mbox{}_{-1}$}
\put(310,25){$E_{k}$}
\put(85,30){\line(0,1){40}}
\put(85,70){\circle*{3}}\put(78,75){$\mbox{}_{-k}$}
\put(88,58){$E_{k+1}$}\put(85,70){\line(1,0){60}}
\put(145,70){\circle*{3}}\put(142,75){$\mbox{}_{-1}$}
\put(142,58){$L_1$}
\put(145,70){\line(1,0){60}}
\put(205,70){\circle*{3}}\put(202,75){$\mbox{}_{-2}$}
\put(202,58){$E_{1}$}
\put(205,70){\line(1,0){15}}\put(230,70){$\dots$}
\put(265,70){\circle*{3}}\put(262,75){$\mbox{}_{-2}$}
\put(260,58){$E_{k-2}$}
\put(265,70){\line(-1,0){15}}\put(265,70){\line(1,0){60}}
\put(325,70){\circle*{3}}\put(320,75){$\mbox{}_{-2}$}
\put(300,58){$E_{k-1}$}\put(325,70){\line(0,-1){40}}
\put(85,30){\line(0,-1){40}}
\put(85,-10){\circle*{3}}\put(80,-17){$\mbox{}_{-2}$}
\put(88,-5){$E_{2k-1}$}
\put(85,-10){\line(1,0){25}}\put(120,-10){$\dots$}
\put(165,-10){\circle*{3}}\put(162,-17){$\mbox{}_{-2}$}
\put(162,-5){$E_{k+2}$}\put(165,-10){\line(-1,0){25}}
\put(165,-10){\line(1,0){80}}
\put(245,-10){\circle*{3}}\put(242,-17){$\mbox{}_{-1}$}
\put(242,-5){$L_{2}$}
\put(245,-10){\line(1,0){80}}
\put(325,-10){\circle*{3}}\put(310,-17){$\mbox{}_{-(k-1)}$}
\put(310,-5){$L_0$}
\put(325,-10){\line(0,1){40}}
\put(200,-33){$\text{Fig.}\, 8$}
\end{picture} \vspace{2cm}

\begin{picture}(300,85)
\put(165,30){\circle*{3}}\put(150,30){$\mbox{}_{-1}$}
\put(168,25){$E_{4}$}
\put(265,30){\circle*{3}}\put(270,30){$\mbox{}_{-1}$}
\put(245,25){$E_{2}$}
\put(265,30){\line(0,1){40}}
\put(165,70){\circle*{3}}\put(158,75){$\mbox{}_{-2}$}
\put(168,58){$E_{3}$}\put(165,70){\line(1,0){50}}
\put(215,70){\circle*{3}}\put(207,75){$\mbox{}_{-1}$}
\put(212,58){$L_1$}
\put(215,70){\line(1,0){50}}
\put(265,70){\circle*{3}}\put(260,75){$\mbox{}_{-2}$}
\put(252,58){$E_{1}$}
\put(265,70){\line(0,-1){80}}
\put(165,70){\line(0,-1){80}}
\put(165,30){\line(0,-1){40}}
\put(165,-10){\circle*{3}}\put(160,-17){$\mbox{}_{-1}$}
\put(168,-5){$L_{2}$}
\put(165,-10){\line(1,0){100}}
\put(265,-10){\circle*{3}}\put(262,-17){$\mbox{}_{-1}$}
\put(252,-5){$L_{0}$}
\put(200,-30){$\text{Fig.}\, 9$}
\end{picture} \vspace{2cm}

Similarly, in  the case when the pencil is given by equation (\ref{P2}), the weighted dual graph  of the curve $C=(\bigcup E_j)\cup(\bigcup L_i)\subset Y$ is depicted in Fig. 10 (remark: {\it if $k=4$, then} $(E_4^2)_Y=-3$).

\begin{picture}(300,85)
\put(85,30){\circle*{3}}\put(70,30){$\mbox{}_{-1}$}
\put(87,25){$E_{2k}$}
\put(325,30){\circle*{3}}\put(330,30){$\mbox{}_{-1}$}
\put(310,25){$E_{k}$}
\put(85,30){\line(0,1){40}}
\put(85,70){\circle*{3}}\put(76,76){$\mbox{}_{-2}$}
\put(87,60){$E_{2k-1}$}\put(85,70){\line(1,0){60}}
\put(124,70){\circle*{3}}\put(109,76){$\mbox{}_{-(k-1)}$}
\put(120,58){$E_{k+1}$}\put(85,70){\line(1,0){60}}
\put(165,70){\circle*{3}}\put(162,75){$\mbox{}_{-1}$}
\put(162,58){$L_1$}
\put(145,70){\line(1,0){60}}
\put(205,70){\circle*{3}}\put(202,75){$\mbox{}_{-2}$}
\put(202,58){$E_{1}$}
\put(205,70){\line(1,0){15}}\put(230,70){$\dots$}
\put(265,70){\circle*{3}}\put(262,75){$\mbox{}_{-2}$}
\put(262,58){$E_{k-3}$}
\put(265,70){\line(-1,0){15}}\put(265,70){\line(1,0){60}}
\put(325,70){\circle*{3}}\put(320,75){$\mbox{}_{-3}$}
\put(300,60){$E_{k-2}$}\put(325,70){\line(0,-1){40}}
\put(85,30){\line(0,-1){40}}
\put(85,-10){\circle*{3}}\put(80,-17){$\mbox{}_{-3}$}
\put(87,-6){$E_{2k-2}$}
\put(127,-10){\circle*{3}}\put(121,-17){$\mbox{}_{-2}$}
\put(119,-6){$E_{2k-3}$}
\put(85,-10){\line(1,0){65}}\put(158,-10){$\dots$}
\put(195,-10){\circle*{3}}\put(192,-17){$\mbox{}_{-2}$}
\put(190,-6){$E_{k+2}$}\put(185,-10){\line(-1,0){10}}
\put(185,-10){\line(1,0){80}}
\put(233,-10){\circle*{3}}\put(230,-17){$\mbox{}_{-1}$}
\put(230,-6){$L_{2}$}
\put(245,-10){\line(1,0){80}}
\put(325,-10){\circle*{3}}\put(320,-17){$\mbox{}_{-2}$}
\put(301,-6){$E_{k-1}$}
\put(275,-10){\circle*{3}}\put(260,-17){$\mbox{}_{-(k-2)}$}
\put(272,-6){$L_0$}
\put(325,-10){\line(0,1){40}}
\put(200,-36){$\text{Fig.}\, 10$}
\end{picture} \vspace{2cm}

Let us prove Theorem 3 in the case when $F\in \mathcal R_{A_{2n-1}}$, $n\geq 2$. The branch curve germ $(B,o)$ of $F$ is given by equation $u^2-v^{2n}=0$, where $o=A_1\in\mathbb C^2\setminus L_2\subset \mathbb P^2$ and $u=\frac{z_0}{z_2}$, $v=\frac{z_1}{z_2}$.
Therefore $\overline B=\overline B_1\cup\overline B_2\subset \mathbb P^2$ is the union of two members (corresponding to $\lambda=1$ and $\mu=\pm 1$) of the pencil given by equation (\ref{P1}).

Let $\tau:Y_1\to Y_2$ be the contraction of the curves $L_1,E_{1},\dots, E_{k-1}$ and $L_2,E_{k+2},\dots$, $E_{2k-1}$ (resp., $L_1,E_1$ and $L_2$ if $n=2$) to  two points.
It is easy to see that $Y_2$ is a smooth surface isomorphic to $E_k\times L_0=\mathbb P^1\times\mathbb P^1$ and $\overline B_1,\overline B_2,E_{k+1}, L_0$ are the fibres of the ruled structure on $Y_2$ defined by the projection to the first factor and $E_k$, $E_{2k}$ are sections.

We have
$$ 
\pi_1:=\pi_1(\mathbb P^2\setminus (\overline B\cup(\cup L_i)))=\pi_1(Y_1\setminus (\overline B\cup(\cup L_i)\cup (\cup E_j)))=
\pi_1(Y_2\setminus (\overline B\cup L_0\cup E_k\cup E_{2k})). $$ 
The natural epimorphism
$$i_*:\pi_1(\mathbb P^2\setminus (\overline B\cup(\cup L_i)))\to
\pi_1(\mathbb P^2\setminus (\overline B\cup L_2))=\pi_1^{loc}(B,o)$$
sends the elements $\gamma_1$ and $\gamma_2$ represented by simple loops around $L_0$ and $L_1$ to the neitral element of $\pi_1^{loc}(B,o)$. Therefore the homomorphism $F_*\circ i_*:\pi_1\to G_F\subset \mathbb S_d$ defines a commutative diagram (\ref{diag2}) of covers in which $\overline F$ is not ramified over $L_0$, $L_1$ and $\overline F_1$ and $\overline F_2$ are not ramified over $L_0$. Now the rest of the proof of Theorem \ref{main3} in the case when $F\in \mathcal R_{A_{2n-1}}$ coincides with the end of the proof in the case when $F\in \mathcal R_{D_4}$.

The proof of Theorem \ref{main3} in the cases when $F\not\in \mathcal R_{A_{2n-1}}\cup\mathcal R_{D_4}$ is similar to one in the case when $F\in \mathcal R_{A_{2n-1}}$. The difference depending on the singularity types of the branch curve germs  $B$ is only in the choice of one of the pencils (\ref{P1}) or (\ref{P2}),  the choice of  $A_1$ or $A_2$ as the point $o$, and the choice of curves contracted to points by $\tau:Y_1\to Y_2\simeq \mathbb P^1\times\mathbb P^1$.

If $F\in\mathcal R_{A_2}$ and the branch curve $B$ is given by $u^2-v^3=0$, , or if $F\in\mathcal R_{D_{5}}$ and the branch curve $B$ is given by $v(u^2-v^{3})=0$, or if $F\in\mathcal R_{E_{7}}$ and the branch curve $B$ is given by $u(u^2-v^{3})=0$, then we use pencil (\ref{P1}) when $n=3$, the point $o$ is $A_2$, and
$u=\frac{z_1}{z_0}$, $v=\frac{z_2}{z_0}$. The morphism $\tau:Y_1\to Y_2$ contracts the curves $L_1\cup E_1\cup E_2$ and $L_0\cup L_2$ to points.

If $F\in\mathcal R_{A_{2n}}$, $n\geq 2$, and the branch curve $B$ is given by $u^2-v^{2n+1}=0$, or if $F\in\mathcal R_{D_{2n+3}}$, $n\geq 2$, and the branch curve $B$ is given by $v(u^2-v^{2n+1})=0$, then we use pencil (\ref{P2}), the point $o$ is $A_1$, and
$u=\frac{z_0}{z_2}$, $v=\frac{z_1}{z_2}$. The morphism $\tau:Y_1\to Y_2$ contracts the curves $L_1\cup E_{k+1}\cup(\bigcup_{j=1}^{k-2}E_j)$ and $L_0\cup L_2\cup (\bigcup_{j=2}^{k-2}E_{k+j})$ to points.

If  $F\in\mathcal R_{D_{2n+2}}$, $n\geq 2$, and the branch curve $B$ is given by $v(u^2-v^{2n})=0$, then we use pencil (\ref{P1}), the point $o$ is $A_1$, and $u=\frac{z_0}{z_2}$, $v=\frac{z_1}{z_2}$. The morphism $\tau:Y_1\to Y_2$ is the same morphism as in the case
$F\in\mathcal R_{A_{2n-1}}$.

If $F\in\mathcal R_{E_6}$ and the branch curve $B$ is given by $u^3-v^4=0$, then we use pencil (\ref{P1}) when $n=4$, the point $o$ is $A_2$, and
$u=\frac{z_1}{z_0}$, $v=\frac{z_2}{z_0}$. The morphism $\tau:Y_1\to Y_2$ contracts the curves $L_1\cup E_1\cup E_2\cup E_3$ and $L_0\cup L_2\cup E_6$ to points.

If $F\in\mathcal R_{E_8}$ and the branch curve $B$ is given by $u^3-v^5=0$, then we use pencil (\ref{P2}) when $n=5$, the point $o$ is $A_1$, and
$u=\frac{z_0}{z_2}$, $v=\frac{z_1}{z_2}$. The morphism $\tau:Y_1\to Y_2$ contracts the curves
$L_1\cup E_1\cup E_2\cup E_3\cup E_6$ and $L_0\cup L_2\cup E_7\cup E_8$ to points.

Now to complete the proof of Theorem \ref{main3}, it suffices to repeat  arguments used in the cases when
$F\in \mathcal R_{A_{2n-1}}\cup\mathcal R_{D_4}$.


\begin{thebibliography}{99}
\bibitem{Ar} V.I. Arnol'd: {\it Normal forms for functions near degenerate critical points, the Weyl groups of $A_k$, $D_k$, $E_k$
 and Lagrangian singularities,} Funct. Anal. Appl., {\bf 6:4} (1972), 254 -- 272.

\bibitem{B} {W. Bart, C. Peters, A. Van de Ven:} {\it Compact complex surfaces,} Springer-Verlag, 1984.

\bibitem{Be} G.V. Belyi: {\it On Galois extensions of a maximal cyclotomic field}, Math. USSR-Izv., {\bf 14:2} (1980), 247 -- 256.


\bibitem{G-R} {H. Grauert, R. Remmert:}
{\it Komplexe R$\ddot{\text{a}}$ume,} Math. Ann., {\bf 136}(1958), 245 -- 318.




\bibitem 
{Ku-3}  {Vik.S. Kulikov:} {\it On the  almost generic covers of the projective plane,} arXiv:1812.01313.

\bibitem{K} {Vik.S. Kulikov:} {\it Dualizing coverings of the plane,} Izv. Math. {\bf 79:5} (2015), 1013 -- 1042.

\bibitem {K-4}  {Vik.S. Kulikov:} {\it On germs of finite morphisms of smooth surfaces,} Proc. Steklov Inst. Math., {\bf 307} (2019).



\bibitem{Mu} D. Mumford:: {\it The topology of normal singularities of an algebraic surface and a criterian for simplisity}, Publ. Math. IHES, no. {\bf 9} (1961).

\bibitem{Sh} I.R. Shafarevich:  {\it Foundations of algebraic geometry}, Izdat. "Nauka'', Moscow, 1972. 567 pp. (in Russian).

\bibitem{Sh1}   I.R. Shafarevich, B.G. Averbukh, Yu.R. Vainberg, A.B. Zhizhchenko, Yu.I. Manin, B.G. Moishezon, G.N. Tyurina, A.N. Tyurin: {\it Algebraic surfaces}, Trudy Mat. Inst. Steklov., 75, Nauka, Moscow, 1965, 3 -- 215.


\bibitem{St}{K. Stein:} {\it Analytische Zerlegungen komplexer R\"aume,}
Math. Ann. \textbf{132} (1956), 63--93.


\bibitem{W} {J.M. Wahl:} {\it Equisingular deformations of plane algebroid curves,} Transactions of the AMS, {\bf 193} (174), 143--170.


\bibitem{Z1} O. Zariski: {\it Studies in singularity. I. Equivalent singularities of plane algebroid curves,}
Amer. J. Math., {\bf 87} (1965), 507--536.


\end{thebibliography}
\end{document}